\documentclass[12pt,a4paper,leqno]{article}

\usepackage{amsfonts}
\usepackage{amsmath}
\usepackage{amssymb}
\usepackage{fullpage}
\usepackage{graphicx}
\usepackage{comment}

\newtheorem{num}{\indent\hskip-4.5pt}[section]
\renewcommand{\thenum}{\rm \arabic{section}.\arabic{num}}
\newcommand{\alku}{\begin{num}\hskip-6pt .\hskip5pt }
\newcommand{\loppu}{\end{num}}

\newcommand{\be}{\addtocounter{num}{1}\begin{equation}}
\newcommand{\ee}{\end{equation}}
\newcommand{\ots}[1]{\bigskip \refstepcounter{num} \indent
\thenum.\hskip 0.2truecm {\bf #1.} \hskip 0.4truecm}

\newcommand{\bea}{\addtocounter{num}{1}\begin{eqnarray}}  
\newcommand{\eea}{\end{eqnarray}}

\font\ff=eusm10 scaled 1200
\font\fff=eusm10 scaled 800
\begin{document}

\def\Im{\mathrm{Im}\,}
\def\Re{\mathrm{Re}\,}
\def\M{\mathcal{M}}
\def\N{\mathbb{N}}
\def\Z{\mathbb{Z}}
\def\R{\mathbb{R}}
\def\C{\mathbb{C}}
\def\B{\mathbb{B}}
\def\L{\mathcal{L}}
\def\Rn{{\mathbb R}^n}
\def\Rns{{\overline{\mathbb R}^n}}
\def\lem{{\bf Lemma.\hskip 0.5truecm}}
\def\pro{{\bf Proposition.\hskip 0.5truecm}}
\def\cor{{\bf Corollary.\hskip 0.5truecm}}
\def\theo{{\bf Theorem.\hskip 0.5truecm}}
\def\rem{{\bf Remark.\hskip 0.5truecm}}
\def\rems{{\bf Remarks.\hskip 0.5truecm}}
\def\exa{{\bf Example.\hskip 0.5truecm}}
\def\defi{{\bf Definition.\hskip 0.5truecm}}
\def\conj{{\bf Conjectures.\hskip 0.5truecm}}
\def\prob{{\bf Open problems.\hskip 0.5truecm}}
\def\E{\hbox{\ff E}}
\def\K{\hbox{\ff K}}
\def\smallE{\hbox{\fff E}}
\def\smallK{\hbox{\fff K}}
\def\mut{\tilde{\mu}}
\def\phit{\tilde{\varphi}}
\def\ACLn{{\mathrm ACL}^n}
\def\proof{{\bf Proof.\hskip 0.5truecm}}       
\def\arth{{\rm arth }}

\newcounter{minutes}\setcounter{minutes}{\time}
\divide\time by 60
\newcounter{hours}\setcounter{hours}{\time}
\multiply\time by 60
\addtocounter{minutes}{-\time}

\begin{center}
{\Large \bf  Generalized Elliptic Integrals and the Legendre $\mathcal{M}$-function}
\end{center}
\medskip

\begin{center}
{\large \sc V. Heikkala, H. Lind\'en,\\ M. K. Vamanamurthy and M. Vuorinen}
\end{center}
\bigskip
\begin{center}
\end{center}

\renewcommand{\thefootnote}{\fnsymbol{footnote}}
\setcounter{footnote}{1}

\renewcommand{\thefootnote}{\arabic{footnote}}
\setcounter{footnote}{0}

\medskip

\begin{abstract}
We study monotonicity and convexity properties of functions arising in 
the theory of elliptic integrals, and in particular in the case of a 
Schwartz-Christoffel conformal mapping from a half-plane to a trapezoid. 
We obtain sharp monotonicity and convexity results for combinations of these 
functions, as well as functional inequalities and a linearization property.  
\end{abstract}

\bigskip


\bigskip


\section{Introduction}

In this paper we continue the study of the modular function 
$\varphi_K^{a,b,c}$ and the generalized modulus $\mu_{a,b,c}$ started in
\cite{hvv}, as well as the generalized elliptic integrals $\K_{a,b,c}$ 
and $\E_{a,b,c}$ (for the notation, see (\ref{eq:muacdef}), 
(\ref{eq:phiackdef}), (\ref{eq:Kdef}) and (\ref{eq:Edef}) below). 
In general, the more freedom the parameter values $a,b$ and 
$c$ are allowed, the more complex and hard-to-handle these functions will be.
As in \cite{hvv} we are here particularly interested in the case 
$b=c-a$. Geometrically this case corresponds to the Schwarz-Christoffel problem
from the unit disk onto a trapezoid, i.e. a quadrilateral with two parallel 
sides (see \cite[Theorem 2.3]{hvv}). In the case $c = 1$, (and $b=1-a$) these 
functions coincide with the special cases $\varphi_K^a$, $\mu_a$, $\K_a$, and 
$\E_a$ which were studied extensively in eg. \cite{AQVV}, and relate to the 
case of a parallelogram.

Given complex numbers
$a,b,$ and $c$ with $c\neq 0,-1,-2, \dots $,
the {\em Gaussian hypergeometric function} is the analytic
continuation to the slit plane $\C \setminus [1,\infty)$ of
the series
\be \label{eq:hypdef}
F(a,b;c;z) = {}_2 F_1(a,b;c;z) = 
\sum_{n=0}^{\infty} \frac{(a,n)(b,n)}{(c,n)} \frac{z^n}{n!}\,,\:\: |z|<1 \,.
\ee
Here $(a,0)=1$ for $a \neq 0$, and $(a,n)$
is the {\em shifted factorial function}
or the {\em Appell symbol}
$$
(a,n) = a(a+1)(a+2) \cdots (a+n-1)
$$
for $n \in \N \setminus \{0\}$, where 
$\N = \{ 0,1,2,\ldots\}$. As usual, we let $\C, \R$ and $\Z$ denote 
respectively, the sets of complex numbers, real numbers, and integers.

A {\em generalized modular equation of order (or degree)} $p>0$ is
\be \label{eq:abc}
\frac{F(a,b;c;1-s^2)}{F(a,b;c;s^2)} = p
\frac{F(a,b;c;1-r^2)}{F(a,b;c;r^2)}\,, \:\: 0 < r < 1\,.
\ee
Sometimes we just call this an $(a,b,c)$-modular equation of order $p$
and we usually
assume that $a,b,c>0$ with $a+b \ge c$, in which case this equation uniquely
defines $s$ as a function of $r$, see \cite[Lemma 4.5]{hvv}. 

Many particular cases of (\ref{eq:abc}) have been studied in the
literature on both
analytic number theory and geometric function theory, 
\cite{BB}, \cite{BBG}, \cite{AVV}, \cite{AQVV}. Rational modular
equations were studied most recently by R. S. Maier in \cite{M}. 
The classical case $(a,b,c)= (\frac{1}{2}, \frac{1}{2}, 1)$
was studied already by Jacobi and many others in the nineteenth century. 
In 1995 B.\ Berndt, S.\ Bhargava, and F.\ Garvan published
an important paper \cite{BBG} in which they studied the case 
$(a,b,c)= (a,1-a,1)$ and $p$ an integer. For several rational
values of $a$ such as $a=\frac{1}{3}, \frac{1}{4}, \frac{1}{6}$ 
and integers $p$ (e.g.\ $p= 2,3,5,7,11,...$)  
they were able to give proofs for numerous algebraic
identities stated by Ramanujan in his unpublished notebooks. 
These identities involve $r$ and $s$ from (\ref{eq:abc}).

To abbreviate (\ref{eq:abc}), we use the decreasing homeomorphism 
$\mu_{a,b,c} : (0,1) \to (0, \infty)$, defined by
\be \label{eq:muacdef}
\mu(r)=\mu_{a,b,c}(r) = \frac{B(a,b)}{2} 
\frac{F(a,b;c;{r'}^2)}{F(a,b;c;r^2)}\,, \:\: r \in (0,1)
\ee
for $a, b, c >0$, $a+b \ge c$, where $B$ is 
the beta function, and $r'$ is the {\em complementary argument}
$r'=\sqrt{1-r^2}$. We call $\mu_{a,b,c}$ the {\em generalized modulus}, 
cf.\ \cite[(2.2)]{LV}.
Now (\ref{eq:abc}) can be rewritten as
\be \label{newmodeq}
\mu_{a,b,c}(s) = p \thinspace \mu_{a,b,c}(r)\,, \:\: 0<r<1\,.
\ee

With $p = 1/K$, $K > 0$, the solution of (\ref{eq:abc}) is then given by
\be \label{eq:phiackdef}
s= \varphi_{K}^{a,b,c}(r) = \mu_{a,b,c}^{-1}(\mu_{a,b,c}(r)/K)\,.
\ee
We call the function $\varphi_K^{a,b,c}$ defined by (\ref{eq:phiackdef})
{\em the $(a,b,c)$-modular function with degree $p= 1/K$}
\cite{BBG}, \cite[(1.5)]{AQVV}. In the case $a<c$ we also use the notation
$$
\mu_{a,c} = \mu_{a,c-a,c}\,, \:\: 
\varphi_{K}^{a,c} = \varphi_{K}^{a,c-a,c}\,.
$$

This article is organized as follows. In Section 2 we introduce the necessary 
notation and the functions studied, as well as known results used in the 
sequel. In Section 3 we obtain various generalizations of
monotonicity results for certain combinations of the generalized elliptic
integrals. The most important results here are Theorems \ref{mutheorem} and
\ref{Mfunctions}, where in particular the latter one concerning the Legendre
$\M$-function leads to many of the results in Section 4. In Section 4 we 
present a number of interesting results, which include the monotonicity 
properties
for functions symmetric with respect to $r$ and $s=\varphi_K^{a,c}(r)$ 
(Lemma  \ref{Ktheo}), the functional inequalities for $\mu_{a,c}$ and
$\varphi_K^{a,c}(r)$, and a linearization result, Theorem \ref{linconj}.
Finally, in Section 5 the dependence on the parameter $c$ for the functions
 $\mu_{a,c}$ and $\varphi_K^{a,c}(r)$ is studied. The main results in this section are Corollary \ref{mudepc} and Theorems \ref{imudpec} and \ref{phidepc}. 
In the final section some open problems are presented.

\section{Preliminaries and definitions}

For $0 < a < \min \{c,1\}$ and $0 < b < c \le a+b $, define the 
{\it generalized complete elliptic integrals 
of the first and second kinds} 
(cf.\ \cite[(1.9), (1.10), (1.3), and (1.5)]{AQVV}) on $[0,1]$ by
\be \label{eq:Kdef}
\K = \K_{a,b,c} = \K_{a,b,c}(r) = \frac{B(a,b)}{2} F(a,b;c;r^2)\,,
\ee
\be \label{eq:Edef}
\E = \E_{a,b,c} = \E_{a,b,c}(r) = \frac{B(a,b)}{2} F(a-1,b;c;r^2)\,,
\ee
\be \label{eq:KpEpdef}
\K' = \K_{a,b,c}' = \K_{a,b,c}(r')\,, \:\: 
\mbox{\rm and} \:\: \E' = \E_{a,b,c}' = \E_{a,b,c}(r')
\ee
for $r \in (0,1)$, $r' = \sqrt{1-r^2}$. The end values are defined by limits
as $r$ tends to $0^+$ and $1^-$, respectively.
In particular, we denote $\K_{a,c} = \K_{a,c-a,c}$ and 
$\E_{a,c} = \E_{a,c-a,c} \,.$ 
Thus, by (\ref{eq:hypas}) below,
$$
\K_{a,b,c}(0) = \E_{a,b,c}(0) = \frac{B(a,b)}{2}
$$
and
$$
\E_{a,b,c}(1) = \frac{1}{2}
\frac{B(a,b) B(c,c+1-a-b)}{B(c+1-a,c-b)}\,,
\:\: \K_{a,b,c}(1) = \infty\,.
$$

Note that the restrictions on $a,b$ and $c$ ensure 
that the function $\K_{a,b,c}$ is increasing and unbounded
whereas $\E_{a,b,c}$ is decreasing and bounded,
as in the classical case $a = b = \frac{1}{2}, c=1$.

Let $\Gamma$ denote Euler's {\em gamma function} 
and let $\Psi$ be its logarithmic derivative (also called the
{\em digamma function}),
$\Psi(z) = \Gamma'(z)/\Gamma(z)$. By \cite[p.\ 198]{Ah} the function
$\Psi$ and its derivative have the series expansions
\be \label{eq:psiser}
\Psi(z) = -\gamma - \frac{1}{z} + \sum_{n=1}^{\infty} \frac{z}{n(n+z)}\,,
\:\: 
\Psi'(z) = \sum_{n=0}^{\infty} \frac{1}{(n+z)^2}\,,
\ee
where $\gamma = -\Psi(1)= \lim_{n \to \infty}(\sum_{k=1}^{n}1/k - \log n) =
0.57721\ldots$ is the {\em Euler-Mascheroni constant}. From
(\ref{eq:psiser}) it is seen that $\Psi$ is strictly increasing on 
$(0, \infty)$ and that $\Psi'$ is strictly decreasing there, so that
$\Psi$ is concave. Moreover,
$\Psi(z+1)= \Psi(z)+ 1/z$ and $\Psi(\frac{1}{2})= -\gamma -2 \log 2$, 
see \cite[Ch. 6]{AS}.

For all $z \in \C \setminus \{0, -1, -2, \ldots\}$ 
and for all $n \in \N$ we have 
\be \label{eq:gamapp}
\Gamma(z+n) = (z,n) \Gamma(z)\,,
\ee
a fact which follows by induction \cite[12.12]{WW}. 
This enables us to extend the Appell symbol for all complex values 
of $a$ and $a+t$, except for non-positive integer values, by
\be \label{eq:appext}
(a,t) = \frac{\Gamma(a+t)}{\Gamma(a)}\,.
\ee
Furthermore, the gamma function satisfies the
reflection formula \cite[12.14]{WW}
\be \label{eq:gamref}
\Gamma(z)\Gamma(1-z) = \frac{\pi}{\sin(\pi z)}
\ee
for all $z \not\in \Z$. In particular, 
$\Gamma(\frac{1}{2})= \sqrt{\pi}$.

The {\em beta function} is defined for $\Re x > 0, \: \Re y > 0$ by
\be \label{eq:betadef}
B(x,y) = \int_{0}^{1} t^{x-1} (1-t)^{y-1} dt
= \frac{\Gamma(x) \Gamma(y)}{\Gamma(x+y)}\,.
\ee
As in this article we are mostly interested in cases where the hypergeometric
parameters satisfy  $0<a<c<1$ and $b=c-a$, we will shorten $B:=B(a,c-a)$ if
no risk for confusion is apparent.

We will make use of the standard notation for contiguous hypergeometric 
functions (cf.\ \cite{R})
$$
F = F(a,b;c;z), \:\: F(a+) = F(a+1,b;c;z), \:\: F(a-) = F(a-1,b;c;z)\,,
$$
etc. We also let 
$$
v = v(z) = F\,, \:\: u = u(z) = F(a-)\,, \:\: v_1 = v_1(z) = v(1-z)\,,
\:\:\mbox{\rm and} \:\: u_1 = u_1(z) = u(1-z)\,.
$$
The behavior of the
hypergeometric function near $z = 1$ in the three cases 
$\Re(a+b-c) < 0$, $a+b = c$, and $\Re(a+b-c) > 0$, 
respectively, is given by
\be \label{eq:hypas}
\left\{ \begin{array}{l}
F(a,b;c;1) = 
\frac{\Gamma(c) \Gamma(c-a-b)}{\Gamma(c-a) \Gamma(c-b)}\,, \\[.5cm]
B(a,b)F(a,b;a+b;z)+\log(1-z) = 
R(a,b)+{\rm O}((1-z)\log(1-z))\,, \\[.4cm]
F(a,b;c;z) = (1-z)^{c-a-b} F(c-a,c-b;c;z)\,,
\end{array} \right.
\ee
where $R(a,b)=-\Psi(a)-\Psi(b)-2 \gamma$.
The above asymptotic formula for the {\it zero-balanced} case $a+b=c$ 
is due to Ramanujan (see \cite{Ask}). This formula is implied by 
\cite[15.3.10]{AS}. Note that $R(\frac{1}{2}, \frac{1}{2}) = \log 16$.

For complex $a,b,c$, and $z$, with $|z| < 1$, we now let
\be \label{eq:Mdef}
\M(z) = \M(a,b,c,z) = 
z(1-z) \left( v_1(z) \frac{dv}{dz} - v(z)\frac{dv_1}{dz} \right) \,.
\ee
Using the Gauss contiguous relations, \cite[p.61]{R}, it is easy to see that
\bea 
\M & = & (c-a)(u v_1 + u_1 v) + (2(a-c)+b)v v_1 \label{eq:Mdef2} \\
 & = & (c-a)(u v_1 + u_1 v - v v_1) + (a+b-c)v v_1 \nonumber
\eea
and that
\bea
(B/2)^2 \M(r^2) &=& (a+b-c)\K \K' + (c-a)[\K \E'+\K'\E - \K \K'].\label{eq:M2def}
\eea

It follows from
\cite[Corollary 3.13(5)]{AQVV} that

\be \label{eq:Maca}
\M(a,1-a,1,r) = \frac{1-a}{\Gamma(a)\Gamma(2-a)} =\frac{ \sin(\pi a) }{ \pi}
\ee
for $0 < a < 1$ and $0 \le r < 1$.
In particular,  we get the classical 
Legendre relation (\cite{AAR}, \cite{BF})
\be \label{eq:Legendre}
\M(1/2,1/2,1,r) = \frac{1}{\pi}\,.
\ee
The function $\M$ will be referred to as {\it the Legendre $\M$-function}, and
it has a central role for the generalizations considered in this 
article. It has the following useful symmetry and convexity properties, 
some of which were established already in \cite[3.17]{hvv} 
(properties (1)-(3)).

\alku \theo \label{Mprop}
For positive constants $a,b,c$ the restriction to $(0,1)$ of 
the continuous function $\M$ has the following properties.\\
(1) $\M(x)=\M(1-x) > 0$ for all  $x \in (0,1)$.\\
(2) If $a+b \le c$, then $\M(x)$ is bounded and extends continuously to
$[0,1]$. In particular, if  $a+b=c=1$, then $\M(x)$ equals the constant
$\sin(\pi a)/\pi$.\\
(3) If $a+b > c$, then $\M$ is unbounded on $(0,1)$ with 
$\M(0^+)=\M(1^-)=\infty$.\\
(4) If $(a+b-1)(c-b)>0,\; a+b \ge c \ge a$ and $ab/(a+b+1)<c$, 
then $\M(a,b,c,r)$ is strictly convex, decreasing in $(0,1/2]$ and increasing 
in $[1/2,1)$.\\
(5) If $(a+b-1)(c-b)<0,\; a+b \le c$, and $ab/(a+b+1) < c$, then 
$\M(a,b,c,r)$ is strictly concave, increasing in $(0,1/2]$ and decreasing 
in $[1/2,1)$.\\
(6) If $a+b \ge c$ then $\M(r) > ab/c$ for all $r \in (0,1)$.
\loppu   

\proof Parts (1)-(3) are proved in the the above mentioned article. 

For (4) and (5) note that by (\ref{eq:Mdef2}) the function $\M$ can be written as
$$\M(a,b,c,r)=(c-a)(u v_1 + u_1 v - v v_1) + (a+b-c)v v_1.$$
In both cases (4) and (5) the constant $(c-a)$ is positive, so 
concavity/convexity of $(c-a)(u v_1 + u_1 v - v v_1)$ follows from the 
assumptions by \cite[2.1]{KV}. The functions $v$ and $v_1$, are both 
log-convex by \cite[1.4]{AVV2}, which follows from the parameter assumption 
$ab/(a+b+1)<c$. Then, so is the product $vv_1$ (by eg. \cite[1.38(5)]{AVV}),
and thus it is convex. Then the convexity/concavity of $(a+b-c)vv_1$ in the
asserted cases also follows. 

For (6), we see that
\begin{eqnarray*}
\M(r) &=& (a+b-c)vv_1+(c-a)[vv_1(a-)+v_1v(a-)-vv_1]\\
&=& (a+b-c)vv_1+ [(c-a)(c-b)/c][(1-r)vv_1(c+)+rv_1v(c+)]\\
&>& (a+b-c)+[(c-a)(c-b)/c] = ab/c.
\end{eqnarray*}
$\square$
\bigskip

Next we record some elementary but useful results for deriving monotonicity 
properties and obtaining inequalities. The first one is the so called 
{\it l'H\^opital's monotone rule}, see \cite[1.25]{AVV} and \cite{AVV3}.


\alku \lem \label{hospital} Let $-\infty < a <b<\infty$, and let $f,g \colon 
[a,b] \to \R$ be continuous on $[a,b]$ and differentiable on $(a,b)$. Let
$g'(x) \ne 0$ on $(a,b)$. Then, if $f'(x)/g'(x)$ is increasing (decreasing)
on $(a,b)$, so are 
$$[f(x)-f(a)]/[g(x)-g(a)]\qquad and \qquad  [f(x)-f(b)]/[g(x)-g(b)].$$
If $f'(x)/g'(x)$ is strictly monotone, then the monotonicity in the conclusion 
is also strict.
\loppu

The second result follows from direct differentiation, and concerns the 
monotonicity  of certain rational functions.


\alku \pro \label{rational}
Assume that $f,g \colon I \to \R$ are differentiable on an interval
$I \subset \R$, and that $a,b,c,d \in \R$. Then 
\begin{eqnarray*}
{\rm sign}\; \left((ad-bc) \frac{d}{dx}\left(\frac{f(x)}{g(x)} \right)\right)
={\rm sign} \frac{d}{dx} \left(\frac{af(x)+bg(x)}{cf(x)+dg(x)}\right).
\end{eqnarray*}
\loppu 

Finally, we record some of the most useful differentiation formulae for the
functions defined in (\ref{eq:muacdef}),(\ref{eq:phiackdef}),(\ref{eq:Kdef}),(\ref{eq:Edef}) and (\ref{eq:Mdef}) (cf. \cite{hvv});

\be \label{eq:Kdiff}
\frac{d \K}{dr} = 
\frac{2}{r {r'}^2} \left((c-a)\E+(b {r}^2+a-c) \K \right)\,,
\ee
\be \label{eq:Ediff}
\frac{d \E}{dr} =
\frac{2(a-1)}{r}\left(\K-\E \right) \,,
\ee
\bea 
\frac{d}{dr}(\K-\E) & = &  
\frac{2}{r {r'}^2} \left(((c-a)+(1-a){r'}^2)\E \right. 
\label{eq:KmEdiff} \\
&& \left. +((a+b) {r}^2-c +{r'}^2) \K \right)\,, \nonumber
\eea
\be \label{eq:KmrpEdiff}
\frac{d}{dr} (\E - {r'}^2 \K) = 
\frac{2}{r}((1-c)\E + (c-1-(b-1)r^2)\K)\,,
\ee
\be \label{eq:muacdiff}
\frac{d}{dr} \mu(r) = 
- \frac{B(a,b)\M(r^2)}{r{r'}^2 v(r^2)^2} =
- \frac{B(a,b)^3\M(r^2)}{4r{r'}^2 \K^2}\,,
\ee
\be \label{eq:phiacdiff}
\frac{\M(s^2)}{\M(r^2)} \frac{ds}{dr} = 
\frac{1}{K} 
\frac{s {s'}^2 v(s^2)^2}{r {r'}^2 v(r^2)^2} = 
\frac{1}{K} 
\frac{s {s'}^2 \K(s)^2}{r {r'}^2 \K(r)^2}\,, 
\;\; s = \varphi_{K}(r)\,,
\ee
\bea \label{eq:Mdiff}
\frac{d\M}{dr} &=& \frac{1}{r(1-r)}\Big((c-a)[(1-c+(a+b-1)r)u(r)v_1(r)\\
&{}&+(-a-b+c+(a+b-1)r)u_1(r)v(r)]\nonumber\\
&{}&+(1-2r)[(c-a)(a+2b-1)-b^2]v(r)v_1(r)\Big).\nonumber
\eea

Note that for the case $(a,b,c)=(1/2,1/2,1)$ the above formulas reduce to the 
classical ones (\cite{BF},\cite{AVV}).

\section{Monotonicity and bounds}

In studying monotonicity and convexity of modular functions, a useful method 
is to combine rational functions consisting of generalized elliptic integrals
whose monotonicity properties are known in different ways. In the following
lemmas we collect some useful properties of such functions, proved in
\cite[4.21, 4.13, 4.24]{hvv}.

\alku \lem \label{EKmonot} For $0< a, b < \min \{c,1\}$ and $c \le a+b \, ,$  
denote $\K = \K_{a,b,c}$
and $\E = \E_{a,b,c}$. Then the function

(1) $f_1(r) = (\K - \E)/(r^2 \K)$ is strictly 
increasing from $(0,1)$ onto $(b/c, 1)$. In
particular, we have the sharp inequality,
$$
\frac{b}{c} < \frac{\K - \E}{r^2 \K} < 1
$$
for all $r \in (0,1)$.

(2) $f_2(r) = (\E-{r'}^2 \K)/r^2$ has positive Maclaurin
coefficients and maps $(0,1)$ onto 
$(B(a,b)(c-b)/(2c), d)$, where
$$ 
 d = \frac{B(a,b)B(c,c+1-a-b)}{2B(c+1-a,c-b)} \, .
$$

(3) $f_5(r) = (r')^{-2} \E$
has positive Maclaurin coefficients and maps
$[0,1)$ onto 
$[B(a,b)/2, \infty)$.  

(4) $f_6(r) = {r'}^2 \K$ 
has negative Maclaurin coefficients, except for the constant term, and maps
$[0,1)$ onto $(0,B(a,b)/2]$.

(5) $f_7(r) = \K$ has positive Maclaurin coefficients and is 
log-convex from $[0,1)$ onto 
$[B(a,b)/2, \infty)$.  In fact, $(d/dr)(\log \K) $ also has
positive Maclaurin coefficients. 

(6) $f_8(r) = (\E - {r'}^2 \K)/(r^2 \K)$ 
is strictly decreasing from $(0,1)$ onto $(0,1-(b/c))$.

(7) $f_9(r) = (\K - \E)/(\E - {r'}^2 \K)$ is
strictly increasing from $(0,1)$ onto $(b/(c-b), \infty)$.
\loppu

\alku \lem \label{EKmonot2} (1) For $0<a<c$ and $b=c-a$, the function 
$h(r) = r^2 \K_{a,c}(r)/\log(1/r')$ is strictly decreasing (respectively, 
increasing)  from $(0,1)$ onto $(1,B(a,b))$ if $a,b \in (0,1)$  (respectively, 
onto $(B(a,b),1)$, if $a, b \in (1,\infty) $).\\
(2) For $0<a,b<c$ and $2ab<c\le a+b < c+1/2$, the function $f(r) = r' \K(r)$ 
is strictly decreasing from $[0,1)$ onto $(0, B(a,b)/2]$.
\loppu

We start with some further monotonicity results for the generalized
elliptic integrals, proved in \cite{AQVV} for the case $c=1,b=1-a$.  
Note that part (1) extends \cite[4.38]{hvv}, as the condition $c \le a+(1/2)$ is not
needed.

\alku \theo \label{hyper} For $c \in (0,1],\; a \in (0,c)$ and $b=c-a$, 
we have that the function\\
(1) $f_1(r)=r \K_{a,c}(r)/\arth(r)$ is strictly decreasing 
from $(0,1)$ onto $(1,B/2)$.\\
(2) $f_2(r)=((B/2)^2-(r'\K_{a,c}(r))^2)/(\E_{a,c}(r)-r'{\thinspace}^2\K_{a,c}(r))$
is strictly increasing from $(0,1)$ onto $(B(c-2ac+2a^2)/(2a),B^2(c-a)/2)$.\\
(3) $f_3(r)=r'{\thinspace}^2(\K_{a,c}(r)-\E_{a,c}(r))/(r^2 \E_{a,c}(r))$ is 
strictly decreasing from $(0,1)$ to $(0,(c-a)/c)$.
\loppu
\medskip

\proof (1) Clearly $f_1(0^+)=B/2$. By l'H\^opital's rule, 
Lemma \ref{hospital}, (\ref{eq:gamapp}), (\ref{eq:betadef}) and the transformation formula and 
evaluation 
at $1$ for hypergeometric functions given in (\ref{eq:hypas}), we see that
\begin{eqnarray*}
f_1(1^-)&=& (B/2) \lim_{r \to 1^-} 2(a/c)(c-a) r'{\thinspace}^2 
F(a+1,c-a+1;c+1;r^2)\\
&=& B (a/c)(c-a) \lim_{r \to 1^-} F(a,c-a;c+1;r^2) =
 B {\textstyle \frac{a}{c}}(c-a) \frac{\Gamma(c+1)}{\Gamma(c-a+1)\Gamma(a+1)}\\
&=&B\ \frac{\Gamma(c)}{\Gamma(c-a)\Gamma(a)} = B \frac{1}{B} =1.
\end{eqnarray*}

Next, let $F_1(r)=rF(a,c-a;c;r^2)$ and $F_2(r)=\arth(r)$. By differentiation we 
get
\begin{eqnarray*}
\frac{F_1'(r)}{F_2'(r)} &=& r'{\thinspace}^2 F(a,c-a;c,r^2)+
2(a/c)(c-a)r^2 F(a,c-a;c+1;r^2)\\
&=&\sum_{n=0}^\infty \frac{(a,n)(c-a,n)}{(c,n)} \frac{r^{2n}}{n!} -
\sum_{n=0}^\infty \frac{(a,n)(c-a,n)}{(c,n)} \frac{r^{2(n+1)}}{n!}\\
&{}& + 2(a/c)(c-a)\sum_{n=0}^\infty \frac{(a,n)(c-a,n)}{(c+1,n)} \frac{r^{2(n+1)}}{n!}\\
&=&1-\sum_{n=1}^\infty \frac{(a,n-1)(c-a,n-1)}{(c,n)n!} 
\cdot [n(1-2a(c-a))-(1-a)(1-c+a)] r^{2n},
\end{eqnarray*}
which is strictly decreasing on $(0,1)$, since
\begin{eqnarray*}
 n(1-2a(c-a))-(1-a)(1-c+a) &\ge& 1-2a(c-a)-(1-a)(1-c+a)\\ 
&>& c-3a(c-a)
 \ge c-{\textstyle \frac{3}{4}} c^2 \ge {\textstyle \frac{1}{4}} c^2>0.
\end{eqnarray*}
Then, by l'H\^opital's rule the function $f$ is also decreasing. 

(2) Let $F(r)=(B/2)^2-(r'\K_{a,c}(r))^2$ and 
$G(r)=\E_{a,c}(r)-r'{\thinspace}^2\K_{a,c}(r)$. Then, using the differentiation
formulas (\ref{eq:Kdiff}) and (\ref{eq:KmEdiff}), we see that
\begin{eqnarray*}
\frac{F'(r)}{G'(r)}&=& \K_{a,c}\; \left(\frac{r^2\K_{a,c}-2(c-a)(\E_{a,c}-r'{\thinspace}^2\K_{a,c})}{ar^2\K_{a,c}+(1-c)(\E_{a,c}-r'{\thinspace}^2\K_{a,c} )}\right).
\end{eqnarray*}
By Lemma \ref{hospital} we need to show that this ratio is strictly increasing. However,
since $\K_{a,c}$ is strictly increasing, and also 
$r \mapsto r^2 \K_{a,c}/(\E_{a,c} -r'{\thinspace}^2 \K_{a,c})$ is, by Lemma \ref{EKmonot}(6), 
the result follows from Proposition \ref{rational} and the fact that 
$(1-c)+2a(c-a) > 0.$ Also, by Lemma \ref{EKmonot}(6)
$$\lim_{r \to 0}  \frac{\E_{a,c} -r'{\thinspace}^2 \K_{a,c}}{r^2 \K_{a,c}} = \frac{a}{c},$$
and so we see that
\begin{eqnarray*}
\lim_{r \to 0^+} \frac{F(r)}{G(r)}&=& \lim_{r \to 0^+} \K_{a,c}\; \left(\frac{1-2(c-a)\frac{\E_{a,c}-r'{\thinspace}^2\K_{a,c}}{r^2\K_{a,c}}}{a+(1-c)\frac{\E_{a,c}-r'{\thinspace}^2\K_{a,c}}{r^2 \K_{a,c}}}\right) =\frac{B}{2} \frac{c-2ac+2a^2}{a}.
\end{eqnarray*}
Furthermore, using the value of $\E_{a,c}(1)$ and the fact that $\lim_{r \to 1^-}
r'\K_{a,c} = 0$, we see that $\lim_{r \to 1^-} F(r)/G(r)=(c-a)B^2/2.$

(3) Follows directly from the fact that $f_3(r)=1-g(r)/\E_{a,c}(r)$, where $g$ is the function $f_2$ in Lemma \ref{EKmonot}(2). $\qquad \square$
\bigskip

The following result extends part of \cite[5.4]{AQVV}.

\alku \lem \label{squareroottimesK}
Let $0<a<c \le 1, B=B(a,b)$ with $b=c-a$, and $\K=\K_{a,c},\; \E=\E_{a,c}$. 
Then the function\\
(1) $f_1(r)=r'{\thinspace}^p\ \K(r)$ is decreasing if and only if $p\ge 2\frac{a}{c}(c-a)$,
in which case $r'{\thinspace}^p\ \K(r)$ is decreasing from $(0,1)$ onto $(0,B/2)$. 
In particular,
$\sqrt{r'} \K(r)$ is decreasing on $[0,1)$.\\
(2) $f_2(r)=r'{\thinspace}^p\ \E(r)$ is increasing  if and only if $p \le 
-\frac{2}{c}(1-a)(c-a)$, in which case it is increasing from $(0,1)$ onto
$(B/2,\infty)$. In particular, $\E(r)/r'{\thinspace}^2$ is increasing on $[0,1)$.
\loppu

\proof (1) Differentiating we get that
$$r(r')^{2-p} f_1'(r)=-pr^2 \K(r)+2(c-a)(\E(r)-r'{\thinspace}^2 \K(r)).$$
This is non-positive if and only if 
$$p \ge 2(c-a) \sup_r \frac{\E(r)-r'{\thinspace}^2 \K(r)}{r^2 \K(r)} =2\frac{a}{c}(c-a),$$
by Lemma \ref{EKmonot}(6). Finally, since $\max\{2\frac{a}{c}(c-a) \  | \ 0<a<c\le1 \} = 1/2$, the function $\sqrt{r'} \K(r)$ will be decreasing for all 
appropriate values of $a$ and $c$. The limiting value at $r=0$ is obvious, and
the one at $r=1$ follows from l'H\^opital's Rule and Lemma \ref{EKmonot}(2).
\medskip

(2) Differentiating yields
$$rf_2'(r)=-p(r')^{p-2} r^2 \E(r)+2(a-1)r'{\thinspace}^p (\K(r)-\E(r)),$$
which is non-negative if and only if
$$-p \ge 2(1-a) \sup_r \frac{r'{\thinspace}^2(\K(r)-\E(r))}{r^2\E(r)} =
\frac{2}{c}(1-a)(c-a),$$
where the value of the supremum follows from Theorem \ref{hyper}(3).
Since $\sup\{\frac{2}{c}(1-a)(c-a) \  | \ 0<a<c\le1 \} = 2$, the function 
$\E(r)/r'{\thinspace}^2$ will be increasing for all 
appropriate values of $a$ and $c$. The limiting values are obvious.
$\square$
\bigskip

\alku \lem \label{logconvexKE}
For $0 < a,b < \min\{c,1\}$ and  $a+b \ge c,\; r \in (0,1)$, we have that
the function\\
(1) $f_1(r)=(r')^{2(a+b-c)} \K_{a,b,c}(r)$ has positive Maclaurin coefficients
and is log-convex on  $(0,1)$ with range $(B(a,b)/2,B(c,a+b-c)/2)$.\\
(2) $f_2(r)=(r')^{2(a+b-c-1)} \E_{a,b,c}(r)$  has positive Maclaurin coefficients
and is log-convex on $(0,1)$ with range $(B(a,b)/2,\infty)$. 
\loppu

\proof (1) From (\ref{eq:hypas}), we have that 
$f_1(r) = (B(a,b)/2) F(c-a,c-b;c;r^2)$, 
so that $(c-a)(c-b) < c(2c-a-b)$ if and only if $ab<c^2+c$, which is true. Hence the 
assertion follows from \cite[ Theorem 3.2(1)]{AVV2}.

(2)  From (\ref{eq:hypas}), we have that $f_2(r) = (B(a,b)/2) F(c+1-a,c-b;c;r^2)$, so that $(c+1-a)(c-b) < c(2c+2-a-b)$, if and only if $(a-1)b<c^2+c$, which 
is true. Hence the assertion follows from \cite[ Theorem 3.2(1)]{AVV2}. 
$\square$
\bigskip

We next derive some monotonicity results for functions combined 
with the $\mu_{a,c}$-function. 

\alku \theo \label{mutheorem}
Let $0<a<c \le 1$. Then the function

(1) $f_1(r) =  \mu_{a,c}(r) + \log r$ is strictly decreasing 
from $(0,1]$ onto $[0, R(a,c-a)/2)$, where $R(a,c-a)$ is as in 
(\ref{eq:hypas}). 

(2) $f_2(r)= \frac{r'{\thinspace}^2 \log r'}{r^2 \log r} \mu_{a,c}(r)$ is 
strictly increasing from $(0,1]$ onto $(1/2,B^2/2]$. 

(3) $f_3(r)= \frac{r' \arth(r)}{r \arth(r')} \mu_{a,c}(r)$ is 
strictly increasing from $(0,1)$ onto $(1,(B/2)^2]$.

(4) $f_4(r)=r'\mu_{a,c}(r)/\log(1/r)$ is strictly increasing from
$(0,1)$ onto $(1,\infty)$. Thus the function 
$\tilde{f}_4(r)=\mu_{a,c}(r)/\log(1/r)$ is also strictly increasing from
$(0,1)$ onto $(1,\infty)$.

(5) $f_5(r)=\mu_{a,c}(r) \arth(r)$ is strictly increasing from
$(0,1)$ onto $(0,(B/2)^2)$.

(6) $f_6(r)= \mu_{a,c}(r) \log (r/r')$ is increasing from
$[1/\sqrt{2},1)$ onto $[0,(B/2)^2)$.
\loppu

\medskip

\proof (1)  Clearly $f_1(1)=0$, and by \cite[1.52(2)]{AVV} it follows that
$f_1(0^+)=R(a,c-a)/2$. From \cite[(4.19)]{hvv} we find that
$$f_1'(r) = \frac{1}{r} -\frac{B(a,c-a)^3 \M(r^2)}{4rr'{\thinspace}^2 
\K_{a,c}^2(r)}= \frac{1}{r}\left(1-\frac{\left(B(a,c-a)/2\right)^2\ B(a,c-a)
 \M(r^2) }{(r' \K_{a,c}(r))^2} \right).$$ 
It now suffices to show that this derivative is negative, which is true if, 
denoting $B=B(a,c-a)$, we have
\bea \label{mulogtoshow}
\frac{(B/2)^2 B \M(r^2)}{(r' \K_{a,c}(r))^2} > 1
\eea
for $r \in (0,1)$. From Lemma \ref{squareroottimesK}(1) it follows that
$g(r)=r'\K_{a,c}(r)$ is strictly decreasing from $[0,1)$ onto
$(0,B/2]$. By Theorem \ref{Mprop} we see that $\M(r^2)$ gets its smallest 
value for $\M(0^+)=\M(1^-) = 1/B$. Then we see that
\begin{eqnarray*}
\frac{(B/2)^2 B \M(r^2)}{(r' \K_{a,c}(r))^2} &\ge& B\; \M(r^2)/r' > B\;
\M(0^+)/r' = \frac{1}{r'} > 1. 
\end{eqnarray*}  
The claim follows.

(2) The function $f_2$ can be rewritten as
$$f_2(r)= \frac{B}{2}\cdot \frac{r'{\thinspace}^2 F(a,c-a;c;r'{\thinspace}^2)}{\log(1/r^2)} 
\cdot \frac{\log(1/r'{\thinspace}^2)}{r^2 F(a,c-a;c;r^2)}.$$
By Lemma \ref{EKmonot2}(1) the second fraction is strictly increasing onto
$(2/B,2]$, and the third onto $(1/2,B/2]$, so the claim follows.

(3) The function $f_3$ can be rewritten as
$$f_3(r)= \frac{B}{2}\cdot \frac{r' F(a,c-a;c;r'{\thinspace}^2)}{\arth(r')} 
\cdot \frac{\arth(r)}{r F(a,c-a;c;r^2)}.$$
Then, in the same way as in part (2), the claim follows from Theorem 
\ref{hyper}(1).

(4) Clearly
$$f_4(r)=\frac{r'F(a,c-a;c;r'{\thinspace}^2) }{\log(1/r)F(a,c-a;c;r^2)} =
\frac{r'{\thinspace}^2F(a,c-a;c;r'{\thinspace}^2)}{\log(1/r)}\cdot \frac{1}{r' F(a,c-a;c;r^2)}.$$
By Lemma \ref{EKmonot2}(1) and part (2) it is then the product of two 
increasing functions. The limiting values also follow immediately. As 
$r \mapsto r'$ is decreasing, $\lim_{r \to 0^+} r' = 1$ and 
$\lim_{r \to  1^-} r' = 0$, the statements for $\tilde{f}_4$ also follow. 

(5) We see that
$$f_5(r)=\frac{B}{2} \frac{\arth(r)}{r \K(r)} r\K'(r).$$
Then it is a product of two increasing functions, by Theorem \ref{hyper}(1)
and Lemma \ref{EKmonot2}(2). Hence, $f_¤$ is increasing itself. The limiting
values are obvious. 

(6) The value $f_6(1/\sqrt{2})=0$ is obvious, while the limit as $r \to 1$ 
follows from (4) and the symmetry property $\mu_{a,c}(r)\mu_{a,c}(r')=
(B/2)^2$. Next, $f_6(r) =(1/2)f_2(r) g(r)$, where 
$$g(r)=\left(\frac{r}{r'}\right)^2 \frac{\log r}{\log r'}\log\left(\frac{r}{r'}\right)^2.$$ 
Hence, by (2) it suffices to prove that
$g(r)$ is increasing on $[1/\sqrt{2},1)$. Put $t=(r/r')^2$, so that 
$$g(r)=\frac{t\ \log t}{\log(t+1)} \log \frac{t+1}{t}$$ 
for $t \in [1,\infty)$. 
Clearly $\log t / \log (t+1)$ is increasing on $[1/\sqrt{2},1)$. Let
$h(t)=t \log ((t+1)/t).$ Then $h'(t)=\log ((t+1)/t)-1/(t+1)$ and 
$h''(t)= -1/(t(t+1)^2) < 0$, so that $h'(t)$ is decreasing. Since
$\lim_{t \to \infty} h'(t)=0$, we get $h'(t) > 0$ on $[1,\infty)$ and thus
$h(t)$ is increasing on $[1,\infty)$. $\qquad \square$
\bigskip

For a quotient of hypergeometric functions with different parameters we obtain
the following results.

\alku \theo \label{differentparams1}
Let $a,b,c,a',b',c'$ be positive constants, satisfying the conditions $a' \ge
a, b' \ge b$, and $c'\le c$, with at least one inequality being strict, and
let $\max\{a',b'\} < c'$. Then the function $f(r) := F(a',b';c';r)/F(a,b;c;r)$ 
is strictly increasing on $[0,1)$ onto $[1,L)$, where 
$$L = \frac{B(c',c'-a'-b') B(c-a,c-b)}{B(c,c-a-b) B(c'-a',c'-b')}$$
in case $a'+b' < c'$, and $L = \infty$ in case $a'+b' \ge c'$.
\loppu

\proof First, $f(0) = 1$ is obvious. Next, let $T_n$ denote the $n$:th
coefficient-quotient, that is $T_n=a_n/b_n$, where $a_n$ and $b_n$ are the $n$:th 
Maclaurin coefficients of $F(a',b';c';r)$ and $F(a,b;c;r)$, respectively.
Then 
$$T_n= \frac{(a',n)(b',n)(c,n)}{(a,n)(b,n)(c',n)},$$
so that $T_{n+1}/T_n = (a'+n)(b'+n)(c+n)/[(a+n)(b+n)(c'+n)] > 1.$ Hence the 
assertion on monotonicity follows from \cite[Theorem 4.3]{hvv}.

Now assume that $a'+b' < c'$. Then $a+b \le a'+b' < c' \le c$, so by
(\ref{eq:hypas}) the assertion for $L$ follows. 

Next, let $a'+b' > c'$, and $a+b>c$. Then $(a'+b'-c')-(a+b-c) =p>0$.
Hence, by (\ref{eq:hypas}), we get 
$$f(r)=(1-r)^{-p} \frac{F(c'-a',c'-b';c';r)}{F(c-a,c-b;c;r)},$$
so that $f(1^-) = \infty$.

Next, if  $a'+b' > c'$ and $a+b=c$, then $a'+b'-c' > a+b-c = 0$, so by 
(\ref{eq:hypas}) it follows that $L=\infty$.

Finally, let $a+b<c$, but $a'+b' \ge c'$. Again, from (\ref{eq:hypas})
it follows that $L=\infty$. $\square$
\bigskip

\alku \cor \label{diffparamscor}
With notation for contiguous hypergeometric functions as in {\rm \cite[p.50]{R}}, 
let $a,b,c$ be positive constants, and let $f = F(a+)/F$, $g = F(b+)/F$ and 
$h = F/F(c+)$. Then $f, g$ and $h$ are all increasing on $[0,1)$, with 
$f(0)=g(0)= h(0) = 1$. Furthermore,\\
(1) $f(1^-) = (c-a-1)/(c-a-b-1)$ if $a+b+1 < c$ and $ = \infty$ otherwise.\\
(2) $g(1^-) = (c-b-1)/(c-a-b-1)$ if $a+b+1 < c$ and $= \infty$ otherwise.\\
(3) $h(1^-) = (c-a)(c-b)/[c(c-a-b)]$ if $a+b < c$ and $= \infty$ otherwise.
\loppu
\bigskip

The particular case $0<a<c<1,\; b=c-a$ requires that we have some knowledge 
about the Legendre $\M$-function, a 
phenomenon which does not show in the case $c=1$, as then $\M(a,c-a,c,r)=\M(a,1-a,1,r)$ is constant by (\ref{eq:Maca}).
In the following theorems we derive some more useful properties of the $\M$-function.

\alku \theo \label{Mextra} Denote $f(r)=(r(1-r))^{a+b-c} \M(a,b,c,r)$. Then the following hold for 
positive $a,b,c$ with $a \le c,\ b \le c$ and $r \in (0,1)$.\\
(1) If $a+b>c$, then the function $f(r)$ is bounded.\\
(2) If $a=c$ or $b=c$, then the function $f(r)$ is the constant $b$ or $a$, 
respectively.\\
(3) If $a+b+1=2c$, then $f(r) =d$, a constant, that is, $\M(r) = 
d (r(1-r))^{1-c}$, where $d=\frac{\Gamma(c)^2}{\Gamma(a)\Gamma(b)}$. 
In particular, $\M(r)$ is constant if and only if $c=1$.
\loppu 

\proof (1) By \cite[3.17(7)]{hvv} we know that the limit of $f$ at $r=0$ is
$(a+b-c)B(c,a+b-c)/B(a,b)$. By symmetry of $f$, it is also the limit at $r=1$. 
Therefore $f$ is bounded if $a+b>c$.

(2) Assume that $c=a$. By (\ref{eq:Mdef2}) we see that
$$f(r)= b (r(1-r))^b\; v(a,b,a,r) v_1(a,b,a,r).$$
By \cite[15.1.8]{AS} we have that
$F(a,b;a,r)=F(b,a;a;r)=(1-r)^{-b}$. Thus
$$f(r)= b\; (r(1-r))^b\ (1-r)^{-b}\ r^{-b} = b,$$
which proves the statement. Since the parameters $a$ and $b$ are interchangeable
in hypergeometric functions, the proof is the same in the case $c=b$.

(3) Let $N(r) = \M(r)/(r(1-r))$. Then, by (\ref{eq:Mdef})
$$N(r) = v_1(r) v'(r) - v(r) v_1'(r),$$ 
and
$$ N'(r) = v_1(r) v''(r) - v(r) v_1''(r).$$ 
As $v$ satisfies the hypergeometric differential equation 
(see \cite[(3),p.54]{R}), we have
$$r(1-r)v''(r) + c(1-2r)v'(r) - abv(r) = 0.$$
Now, $v_1'(r) = - v'(1-r)$ and $v_1''(r) = v''(1-r)$. 
Hence,
$$r(1-r) v_1''(r) +c (1-2r) v_1'(r)- ab v_1(r) = 0,$$ 
and thus
$$r(1-r) N'(r) + c(1-2r) N(r) = 0.$$
Hence
$$\frac{d}{dr}\left[ (r(1-r))^{c} N(r)\right] = 0 = \frac{d}{dr}\left[ 
(r(1-r))^{c-1} \M(r)\right],$$ 
so that $\M(r) = d\ (r(1-r))^{1-c}$, where $d$ is a constant. 

We now show that $d = \Gamma(c)^2/(\Gamma(a)\Gamma(b))$. Taking the limit 
as $r \to 0^+$, we have $d = f(0^+)$.
{\it Case (i):} $c = 1$, so that $b = 1-a$. Then by (\ref{eq:hypas})
$$f(r) = r(1-r)v_1(r)v'(r) + a(1-a)(1-r)v(r)v(1-a,a;2;1-r),$$ 
so that 
\begin{eqnarray*}
d&=& a(1-a)v(1-a, a;2;1)= a(1-a) \frac{\Gamma(2)\Gamma(1)}{\Gamma(1+a)\Gamma(2-a)}
= \frac{1}{\Gamma(a)\Gamma(1-a)},
\end{eqnarray*}
as required. Note that in this case $d= \sin(\pi a)/ \pi$.\\
{\it Case (ii):} $0 < c < 1$. In this case we have $0 < a+b < c < 1$. 
Then
\begin{eqnarray*}
f(r) &=& (r(1-r))^{c} \left(v_1(r) v'(r) + (ab/c)v(r)v(a+1,b+1;c+1;1-r)\right)\\
&=& (r(1-r))^{c} \left[v_1(r) v'(r) + (ab/c)r^{-c} v(r)v(c-a,c-b;c+1;1-r)\right],
\end{eqnarray*}
so that 
$$f(0^+) = 0 + (ab/c) v(c-a,c-b:c+1;1) = (ab/c) \frac{\Gamma(c+1)\Gamma(c)}{\Gamma(a+1)\Gamma(b+1)}
= \frac{\Gamma(c)^2}{\Gamma(a)\Gamma(b)}.$$
{\it Case (iii):} $c > 1.$ This is similar to case (ii). $\square$
\bigskip

The following corollary is a direct consequence of Theorem \ref{Mextra}(3) and
the formulas (\ref{eq:muacdiff}) and (\ref{eq:phiacdiff}).

\alku \cor \label{Mcorollary}  Let $\mu = \mu_{a,b,c}$ and let 
$s = \varphi_K^{a,b,c}(r)$. If $a,b,c$ are positive with $a \le c$ and
$b\le c,\; r \in (0,1)$ and $a+b+1=2c$, then
we have the following generalized derivative formulas.\\
(1) $$\frac{d\mu}{dr} = - \frac{D}{r^{2c-1} r'^{2c}\K(r)^2},$$ 
where $D= \frac{(\Gamma(a)\Gamma(b)\Gamma(c))^2}{4\Gamma(a+b)^3}$.\\
(2) $$\frac{ds}{dr} = \frac{1}{K} \left(\frac{s}{r}\right)^{2c-1}
\left(\frac{s'}{r'}\right)^{2c} \left(\frac{\K(s)}{\K(r)} \right)^2.$$\\
\loppu
\bigskip

\alku \theo \label{Mfunctions} Let $0<a<c\le 1,\; b=c-a$ and $\M(r)=\M(a,c-a,c,r)$. Then\\
(1) The inequality
$$\M(r^2)-2r^2\M'(r^2) \ge (c-a)a > 0$$
holds for all $r \in [0,1]$. In particular the function $f(r)=r/\M(r^2)-a(c-a)r$ is increasing from $[0,1]$ onto
$[0,B-a(c-a)]$.\\
(2) The function $g(r)=f(r')$ is decreasing from $[0,1]$ onto $[0,B-a(c-a)]$.
\loppu

\proof (1) First, if $c=1$, then $\M(r^2)$ is a positive constant, hence the
assertion is trivial. We then assume that $0<c<1$. In this case, $(a+b-1)(c-b)=
(c-1)a < 0$, so by Theorem \ref{Mprop} $\M'(r^2)>0$ for $r \in (0,1/\sqrt{2})$
and  $<0$ for $r \in (1/\sqrt{2},1)$. Let
$$F_1=F(c;r'{\thinspace}^2)F(c+;r^2) \qquad {\rm and}\qquad F_2=F(c;r^2)F(c+;r'{\thinspace}^2),$$
where the parameter triple of $F$ is $(a,c-a;c)$. Then we see that both $F_1$ 
and $F_2$ are nonnegative, and in fact $\ge 1$.
As in \cite{KV} (11) and (27), we see that
\begin{eqnarray*}
\M(r^2)= (c-a)\frac{a}{c}\Big(r^2 F_1+r'{\thinspace}^2 F_2\Big)\quad {\rm and} \quad
\M'(r^2)=(c-a)\frac{a}{c}(1-c)\Big(F_1- F_2\Big).
\end{eqnarray*}
Now $\M'(r^2)$ is negative in $(1/\sqrt{2},1)$, so from the equation above we see that in this interval
$F_1-F_2$ is also negative. Then
\begin{eqnarray*}
\M(r^2)^2 \frac{d}{dr} \frac{r}{\M(r^2)} &=&\M(r^2)-2r^2\M'(r^2) = 
(c-a)\frac{a}{c}\Big(r^2F_1+r'{\thinspace}^2F_2-2(1-c)r^2(F_1-F_2)\Big)\\
&\ge&(c-a)\frac{a}{c}\Big(r^2F_1+r'{\thinspace}^2F_2\Big) \ge (c-a)\frac{a}{c}(r^2+r'{\thinspace}^2)= \frac{(c-a)a}{c}.
\end{eqnarray*}
In the case $r \in (0,1/\sqrt{2})$ and  $c \ge 1/2$ both $F_1-F_2$ and 
$(2c-1)$ are nonnegative. Then we see that
\begin{eqnarray*}
\M(r^2)-2r^2\M'(r^2) &=& (c-a)\frac{a}{c}\Big(r^2F_1+r'{\thinspace}^2F_2-2(1-c)r^2(F_1-F_2)\Big)\\
&=&(c-a)\frac{a}{c}\Big(F_2 + (2c-1)(F_1 - F_2)r^2]\Big)\\
&\ge& (c-a)\frac{a}{c} F_2 \ge  (c-a)\frac{a}{c}.
\end{eqnarray*}
For $c \le 1/2$ the
expression $(2c-1)$ is non-positive. Thus, using (\ref{eq:hypas}), and
the inequality $r F(a,c-a;c;r'{\thinspace}^2) \le 1$ which follows from Lemma 
\ref{EKmonot2}(2), we get
\begin{eqnarray*}
&{}&(c-a)\frac{a}{c}\Big(r^2(2c-1)F_1+(1-r^2(2c-1))F_2\Big)\\
&\ge&(c-a)\frac{a}{c}\Big((2c-1)\frac{c}{(c-a)a}\frac{1}{B}r+1-r^2(2c-1)\Big)\\
&\ge&(2c-1)\frac{1}{B}r+r^2 \frac{a(c-a)(1-2c)}{c}+\frac{(c-a)a}{c}.
\end{eqnarray*}
Finally, using the inequality $1/B(a,b) \le 2ab/(a+b)$ (see \cite[1.50]{AVV}) 
together with the fact that $r(1-r) \le 1/4$, we obtain
\begin{eqnarray*}
\M(r^2)^2 \frac{d}{dr}\frac{r}{\M(r^2)} &\ge& \frac{(c-a)a}{c}\Big(1-2(1-2c)r+(1-2c)r^2 \Big)\\
&=&\frac{(c-a)a}{c} \Big(2c+(1-2c)-2(1-2c)r+(1-2c)r^2 \Big)\\
&=&\frac{(c-a)a}{c} \Big(2c+(1-2c)(1-r)^2\Big) \ge 
\frac{(c-a)a}{c} 2c=2a(c-a). 
\end{eqnarray*}
This proves the statement.

Part (2) follows directly from the equality $\M(x)=\M(1-x)$ by interchanging $x$ with $x'$ in 
part (1). $\square$
\bigskip
\bigskip

\section{Functional inequalities and linearization}
\label{sect:linearization}

In this section we generalize the functional inequalities for the modular 
function $\varphi_K^a(r)$ proved in \cite{AQVV} to hold also for the 
generalized modular function $\varphi_K^{a,b,c}(r)$ in the case $b=c-a$.
We start by a generalization of the results in \cite[6.2]{AQVV}.

\alku \lem \label{Ktheo} Let $a<c \le 1,\; K \in (1,\infty),\; r \in (0,1)$,
and let $s=\varphi_K^{a,c}(r)$ and $t=\varphi_{1/K}^{a,c}(r)$. Then the function\\
(1) $f_1(r)=s/r$ is decreasing from $(0,1)$ onto $(1,\infty)$,\\
(2) $f_2(r)=s'/r'$ is decreasing from $(0,1)$ onto $(0,1)$,\\
(3) $f_3(r)=\K(s)/\K(r)$ is increasing from $(0,1)$ onto $(1,K)$,\\
(4) $f_4(r)=\K'(s)/\K'(r)$ is increasing from $(0,1)$ onto $(1/K,1)$,\\
(5) $f_5(r)=s'\K_{a,c}(s)^2/(r'\K_{a,c}(r)^2)$ is decreasing from $(0,1)$ onto $(0,1)$,\\
(6) $f_6(r)=s\K'_{a,c}(s)^2/(r\K'_{a,c}(r)^2)$ is decreasing from $(0,1)$ onto $(1,\infty)$,\\
(7) $g_1(r)=t/r$ is increasing from $(0,1)$ onto $(0,1)$,\\ 
(8) $g_2(r)=t'/r'$ is increasing from $(0,1)$ onto $(1,\infty)$,\\
(9) $g_3(r)=\K(t)/\K(r)$ is decreasing from $(0,1)$ onto $(1/K,1)$,\\
(10) $g_4(r)=\K'(t)/\K'(r)$ is decreasing from $(0,1)$ onto $(1,K)$,\\
(11) $g_5(r)=t'\K_{a,c}(t)^2/(r'\K_{a,c}(r)^2)$ is increasing from $(0,1)$ onto $(1,\infty)$,\\
(12) $g_6(r)=t\K'_{a,c}(t)^2/(r\K'_{a,c}(r)^2)$ is increasing from $(0,1)$ onto $(0,1)$.\\
\loppu

\proof (1) Differentiating we see that
\begin{eqnarray*}
f_1'(r) = \frac{ss'{\thinspace}^2 \K(s)\K'(s)\M(r^2)}{rr'{\thinspace}^2 \K(r)\K'(r)\M(s^2)}\cdot r -s \le 0
\end{eqnarray*}
if and only if 
\bea \label{help2}
 \frac{s'{\thinspace}^2 \K(s)\K'(s)}{\M(s^2)} \le  \frac{r'{\thinspace}^2 \K(r)\K'(r)}{\M(r^2)}.
\eea
As $\mu_{a,c}(s) = \mu_{a,c}(r)/K$, we see that $s > r$ for all $r \in (0,1)$, 
and thus (\ref{help2}) holds if $x \mapsto x'{\thinspace}^2 \K(x)\K'(x)/\M(x^2)$ is 
decreasing, which is true by Theorem \ref{Mfunctions}(2) and Lemma 
\ref{EKmonot2}(2). The limiting value at $1$ is clear. For the limiting value 
at $0$ we see that since $s/r = (s/r^{1/K})(1/r^{1-1/K})$, we get
\begin{eqnarray*}
\log(s/r)&=&(1-1/K)\log(1/r)+(\log s - (1/K) \log r)\\
&=& (1-1/K)\log(1/r)+\big((\mu(s)+\log s)-(1/K)(\mu r +\log r)\big),
\end{eqnarray*} 
which by Theorem \ref{mutheorem}(1) tends to $\infty$.
\medskip

(3) Differentiating we have that
\begin{eqnarray*}
\K_{a,c}(r)^2 f_1'(r) &=&2(c-a)\left[\K_{a,c}(r) \frac{\E_{a,c}(s)-s'{\thinspace}^2 \K_{a,c}(s)}{ss'{\thinspace}^2} \frac{ds}{dr}-\K_{a,c}(s)\frac{\E_{a,c}(r)-r'{\thinspace}^2 \K_{a,c}(r)}{rr'{\thinspace}^2}  \right]\\
&=&\frac{2(c-a)\K_{a,c}(s)\M(r^2)}{rr'{\thinspace}^2\K'_{a,c}(r)}
\bigg[\frac{\K'_{a,c}(s)(\E_{a,c}(s)-s'{\thinspace}^2 \K_{a,c}(s))}{\M(s^2)}\\
&{}&-\frac{\K'_{a,c}(r)(\E_{a,c}(r)-r'{\thinspace}^2 \K_{a,c}(r))}{\M(r^2)} \bigg].
\end{eqnarray*}
Then, since by Theorem \ref{phiperr} $r< r^{1/K}< s$, it suffices to show that
\begin{eqnarray*}
\frac{\K'_{a,c}(x)(\E_{a,c}(x)-x'^2 \K_{a,c}(x))}{\M(x^2)}=
\frac{x^2 \K'_{a,c}(x)}{\M(x^2)}\cdot \frac{\E_{a,c}(x)-x'^2 \K_{a,c}(x)}{x^2}.
\end{eqnarray*}
is increasing. But this follows from  Theorem \ref{Mfunctions}(1) together with
parts (2) and (9) of Theorem \ref{EKmonot2}(2). The limiting values are clear. 
\medskip

(5) Differentiating, we see that
$f_3'(r)$ is negative if and only if the function
$$F(x)=\frac{x^2 \K_{a,c}(x)\K'_{a,c}(x)}{\M(x^2)} 
\left(1-4(c-a)\frac{\E_{a,c}(x)-x'^2 \K_{a,c}(x)}{x^2\K_{a,c}(x)} \right)$$
is increasing. But this follows from Lemma \ref{EKmonot}(6) 
together with \ref{Mfunctions}(3). The limiting values follow from
the limiting values in parts (2) and (3), as
$$\lim_{r \to 0} \frac{s'\K(s)^2}{r'\K(r)^2} = \lim_{r \to 0}\frac{s'}{r'} 
\cdot \lim_{r \to 0} \frac{\K(s)^2}{\K(r)^2} = 1 \cdot 1 =1.$$
and
$$\lim_{r \to 1} \frac{s'\K(s)^2}{r'\K(r)^2} = \lim_{r \to 1}\frac{s'}{r'} 
\cdot \lim_{r \to 1} \frac{\K(s)^2}{\K(r)^2} = 0 \cdot K^2 =0.$$
\medskip

As $f_2(r)=1/f_1(s'),\; f_4(r)=1/f_3(s')$ and $f_6(r)=1/f_5(s')$, parts (2), (4) and (6) follow. The parts (7)-(12) follow from (1)-(6), as $g_i(r)=1/f_i(t)$ 
for $i=1,2,3,4,5,6$. $\square$
\bigskip

We continue by proving some functional inequalities for the function $\mu_{a,c}$. 

\alku \theo \label{mufunc} Let $0<a<c \le 1$. Then, denoting 
$f(r)=\mu_{a,c}(r)=\mu(r)$, the
function $g_1(r)=(1-r)f'(r)$ is increasing, and the function $g_2(r)=rf'(r)$ is 
decreasing. In particular, the inequalities
\begin{eqnarray*}
\mu_{a,c}\left(1-\sqrt{(1-u)(1-t)}\right) &\le& \frac{\mu_{a,c}(u)+\mu_{a,c}(t)}{2} \le \mu_{a,c}(\sqrt{ut})
\end{eqnarray*} 
hold for all $u,t \in (0,1)$ with equality if and only if $u=t$.
\loppu

\proof We first see that for the function $g_1(r)$
\begin{eqnarray*}
-g_1(r)&=&\frac{B^3}{4}\frac{\M(r^2)}{r} \frac{1}{(1+r)\K(r)^2}.
\end{eqnarray*} 
Clearly this is decreasing by Theorem \ref{Mfunctions}(1), so 
that $g_1(r)$ is increasing. Also
\begin{eqnarray*}
-g_2(r)&=&\frac{B^3}{4}\frac{\M(r^2)}{r'} \frac{1}{r'\K(r)^2},
\end{eqnarray*} 
which is increasing by Theorems \ref{Mfunctions}(2) and 
\ref{squareroottimesK}(1), so that $g_2(r)$ is decreasing. These
monotone properties imply that that the function $f(1-e^{-t})$ is 
convex on $(0,\infty)$ and that the function $f(e^{-t})$ is concave
on $(0,\infty)$, and so the asserted inequalities follow.$\qquad \square$
\bigskip

\alku \theo \label{phiperr} For each $0<a<c \le 1$ and $K > 1,$ the  function
$f(r) = \varphi^{a,c}_K(r) / r^{1/K}$ is strictly decreasing from 
$(0,1]$ onto $[1, e^{(1-(1/K))R(a,c-a)/2}).$ In particular,
$$r^{1/K} < \varphi^{a,c}_K(r) < e^{(1-(1/K))R(a,c-a)/2}  r^{1/K}.$$
Also, the function $g(r) = \varphi^{a,c}_{1/K}(r) / r^{K}$ is strictly increasing 
from $(0,1]$ onto $(e^{(1-K)R(a,c-a)/2},1].$ In particular
$$r^{K} > \varphi^{a,c}_{1/K}(r) > e^{(1-K)R(a,c-a)/2}  r^{K}.$$
\loppu

\proof If $s = \varphi^{a,c}_K(r),$ then $\mu_{a,c}(s) =  \mu_{a,c}(r)/K,$ 
and $s > r,$ for all $r \in (0,1)$ and $K > 1$. Differentiating we get
\begin{eqnarray*}
\frac{f'(r)}{f(r)} =\frac{1}{Kr}\left(\left(\frac{s'\K(s)}{r'\K(r)}\right)^2 \frac{\M(r^2)}{\M(s^2)}-1 \right).
\end{eqnarray*}
This derivative is negative if and only if $(s'{\thinspace}^2\K(s)^2)/\M(s^2) \le  (r'{\thinspace}^2\K(r)^2)/\M(r^2)$, 
that is, if the function $x \mapsto (x'{\thinspace}^2 \K(x)^2)/\M(x^2)$ is decreasing. This, however, follows
from Theorems \ref{Mfunctions}(1) and \ref{squareroottimesK}(1), as 
$$\frac{x'^2\K(x)^2}{\M(x^2)} = \frac{x'}{\M(x^2)} (\sqrt{x'} \K(x))^2.$$
Then $f$ is indeed strictly decreasing. By Theorem \ref{mutheorem}(1) 
\begin{eqnarray*}
\log(s/r^{1/K}) &=&[\mu(s)+\log(s)]-(1/K)[\mu(r)+\log(r)]
\end{eqnarray*}
tends to $(1-(1/K)) R(a,c-a)/2$, as $r \to 0$. The proof for the function 
$g$ follows the same pattern.$\qquad \square$
\bigskip

\alku \rem 
{\rm 
We observe that in Theorem \ref{phiperr} for $a = 1/2,$ and $c = 1,$ the 
coefficient in the upper bound reduces to the classical constant $4^{1-(1/K)}$ 
\cite{LV}.
}
\loppu

\alku \theo \label{funcineq1} Let $0<a<c \le 1$ and $K \in (1,\infty)$. Then
the function\\
(1) the function $f_1(r)=\log(\varphi_K(r'))$ is decreasing and 
concave on $(0,1)$. In particular
$$\varphi_K(u')\varphi_K(t') \le \varphi_K\left(\sqrt{1-\left(\frac{u+t}{2} \right)^2}\right)^2,$$
and
$$\varphi_K(u)\varphi_K(t) \le \varphi_K\left(\sqrt{1-\sqrt{(1-u^2)(1-t^2)}}\right)^2.$$
for all $u,t \in (0,1)$, with equality if and only if $u=t$.\\
(2) The function $f_2(r)=\log(\varphi_K(r'{\thinspace}^2))$ is decreasing and 
concave 
on $(0,1)$. In particular
$$\varphi_K(u'{\thinspace}^2)\varphi_K(t'{\thinspace}^2) \le \varphi_K\left(1-\left(\frac{u+t}{2} \right)^2\right)^2,$$
and
$$\varphi_K(u)\varphi_K(t) \le \varphi_K\left(1-\sqrt{(1-u^2)(1-t^2)}\right)^2.$$
for all $u,t \in (0,1)$, with equality if and only if $u=t$.\\
(3) The function $f_3(r)=\log(\varphi_K(1-e^{-r}))$ is increasing and concave 
on $(0,\infty)$.
In particular
$$\varphi_K(1-u)\varphi_K(1-t) \le \varphi_K(1-\sqrt{ut})^2,$$
and
$$\varphi_K(u)\varphi_K(t) \le \varphi_K\left(1-\sqrt{(1-u)(1-t)}\right)^2.$$
for all $u,t \in (0,1)$, with equality if and only if $u=t$.
\loppu

\proof (1) Denote $t=\varphi_{1/K}(r)$. Then we see that 
$\varphi_K(r')=\sqrt{1-\varphi_{1/K}(r)^2} = \sqrt{1-t^2} = t'$. Now
\begin{eqnarray*}
\frac{d(t')}{dr}&=& \frac{tt'{\thinspace}^2\K(t)\K'(t)\M(r^2)}{rr'{\thinspace}^2\K(r)\K'(r)\M(t^2)}
\cdot \left(-\frac{t}{t'} \right) = -\frac{t't^2\K(t)\K'(t)\M(r^2)}
{rr'{\thinspace}^2\K(r)\K'(r)\M(t^2)}
= -\frac{1}{K} \frac{t't^2 \K'(t)^2 \M(r^2)}{rr'{\thinspace}^2\K'(r)^2 \M(t^2)}.
\end{eqnarray*}
Thus
\begin{eqnarray*}
\frac{df_1}{dr}&=& \frac{1}{t'}\frac{d(t')}{dr} = -\frac{1}{K} \left(\frac{t\K'(t)^2}{r\K'(r)^2}\right)\left(\frac{t}{\M(t^2)}\right)\left(\frac{\M(r^2)}{r'{\thinspace}^2}\right),
\end{eqnarray*}
where each of the bracketed functions is positive and increasing, by  
Theorems \ref{Ktheo}(12) and \ref{Mfunctions}. Thus $df_1/dr$ is 
negative and decreasing, and $f_1$ is decreasing and concave. Then the 
convexity inequality $f((x+y)/2) \ge (f(x)+f(y))/2$ directly yields the first
inequality. The rewritten inequality follows from change of variables.

(2) Again, let $t=\varphi_{1/K}(r)$, and $u(r)=\sqrt{2r^2-r^4}$. The function
$u(r)$ is easily shown to be increasing in $(0,1)$. Now
\begin{eqnarray*}
\frac{df_2}{dr}&=& \frac{t'(u)t(u)^2\K(t(u))\K'(t(u))\M(u^2)}{u'u^2\K(u)\K'(u)\M(t(u)^2)} \cdot (-2r) \cdot \frac{1}{t'(u)}\\
&=& -\frac{1}{K} \left(\frac{t(u)\K'(t(u))^2}{u\K'(u)^2}\right)\left(\frac{t(u)}{\M(t(u)^2)}\right)\left(\frac{\M(u^2)}{u'}\right)\left(\frac{2r}{\sqrt{2r^2-r^4}}\right).
\end{eqnarray*} 
Also here all the bracketed functions are positive and increasing, and thus 
$df_2/dr$ is negative and decreasing, and $f_1$ is decreasing and concave.
The rest of the statement is proved as in (1).

(3) With $x=1-e^{-r}$ and $s=\varphi_K(x)$ we have
$$f_3'(r)=\left(\frac{1-x}{Kx}\right) \left(\frac{s\K(s)^2}{x\K(x)^2}\right) 
\left(\frac{\M(x^2)}{x'}\right)\left(\frac{s}{\M(s^2)}\right),$$
which is decreasing by Theorem \ref{Ktheo}(6). The rest of the statement 
is proved as in the previous cases.$\qquad \square$

\bigskip

\alku \theo \label{linconj} Let $p:(0,1) \to (-\infty, \infty)$ and $q:(-\infty, \infty) \to (0,1)$
be given by $p(x) = 2 \log(x/x')$ and $q(x) = p^{-1}(x) = 
\sqrt{e^x/(e^x + 1)}$, respectively, and for $a \in (0,1)$, $c \in (a,1]$,
$K \in (1, \infty)$, let $g, h:(-\infty, \infty) \to (-\infty, \infty)$
be defined by $g(x) = p(\varphi^{a,c}_K(q(x)))$ and
$h(x) = p(\varphi^{a,c}_{1/K}(q(x)))$. Then
$$
g(x) \geq \left\{ \begin{array}{ll}
                  Kx\,, & \mbox{\it if $x \geq 0$} \\
                  \frac{x}{K}\,, & \mbox{\it if $x < 0$}
                  \end{array}
          \right.
\:\:\:\: \mbox{\it and} \:\:\:\:
h(x) \leq \left\{ \begin{array}{ll}
                  \frac{x}{K}\,, & \mbox{\it if $x \geq 0$} \\
                  Kx\,, & \mbox{\it if $x < 0$.}
                  \end{array}
          \right.
$$
\loppu

\medskip

\proof First, if  $x>0$, then
\begin{eqnarray*}
g(x) \ge Kx &\Leftrightarrow & \varphi_K^{a,c}(q(x)) \ge q(Kx) \\[.3cm]
&\Leftrightarrow & \mu_{a,c}^{-1}\left(\frac{1}{K} \mu_{a,c} 
\left(\sqrt{\frac{e^x}
{e^x+1}}\right) \right)
\ge \sqrt{\frac{e^{Kx}}{e^{Kx}+1}} \\[.2cm]
&\Leftrightarrow & \mu_{a,c}\left(\sqrt{\frac{e^x}{e^x+1}}\right)
\le K\mu_{a,c}\left(\sqrt{\frac{e^{Kx}}{e^{Kx}+1}}\right)\,.
\end{eqnarray*}
This will be true if $f(K) =
K\mu_{a,c} (\sqrt{e^{Kx}/(e^{Kx}+1})$ is increasing on $[1, \infty)$.

Now, setting $r = \sqrt{e^{Kx}/(e^{Kx}+1)}$, we have
$r^2= e^{Kx}/(e^{Kx}+1)$, and $r'{\thinspace}^2 = 1/(e^{Kx}+1)$.
Then $f(K)=(2/x) f_6(r)$, where $f_6$ is as in Theorem \ref{mutheorem}(6),
and thus increasing, as $r(K)$ is increasing as a function of $K$.

Let still $x > 0$. Then
\begin{eqnarray*}
g(-x) \geq -x/K
&\Leftrightarrow & \varphi_K^{a,c} \left(\sqrt{\frac{e^{-x}}{e^{-x}+1}}\right) \ge
\sqrt{\frac{e^{-x/K}}{e^{-x/K}+1}} \\[.2cm]
& \Leftrightarrow &
\mu_{a,c}^{-1}\left(\frac{1}{K}\mu_{a,c}
\left(\sqrt{\frac{1}{e^x+1}}\right)\right) \ge
\sqrt{\frac{1}{e^{x/K}+1}} \\[.2cm]
& \Leftrightarrow  &\mu_{a,c}\left(\frac{1}{\sqrt{e^x+1}}\right) \le
K \mu_{a,c} \left(\frac{1}{\sqrt{e^{x/K}+1}}\right)\,.
\end{eqnarray*}
This is true if $F(K) = K\mu_{a,c}(1/\sqrt{e^{x/K}+1})$ is increasing
on $[1, \infty)$. Let $t = 1/\sqrt{e^{x/K}+1}$. Then $t \in (0, 1/\sqrt{2})$
and $t^2 = 1/(e^{x/K}+1)$, $t'{\thinspace}^2 = e^{x/K}/(e^{x/K}+1)$, 
$x = 2K \log(t'/t)$.
Now $f(K)=(B^2/8)(x/f_6(t'))$, where $f_6$ is as in Theorem 
\ref{mutheorem}(6), and thus increasing, as $t'(K)$ is decreasing as a 
function of $K$. Finally, the proof of $h(x)$ is similar.$\qquad \square$
\bigskip
\bigskip

\section{Dependence on $c$}
\label{sect:depc}

In this section we study how the functions $\mu_{a,c}, \mu^{-1}_{a,c}$ and 
$\varphi_K^{a,c}$ depend on the parameter $c$. corresponding results for 
the case $c=1$ can be found in the articles \cite{AQVV} and \cite{QV1}.

\ots{Notation} \label{pabnot}
For $0 < a < c$ and 
$t > 0$ we denote
$$
P(a,c,t) = \Psi(c-a+t) - \Psi(c+t)\,,
$$
$$
A = A_t = A(a,c,t) = \frac{(c-a,t)}{(c,t)} = 
\frac{\Gamma(c-a+t)\Gamma(c)}{\Gamma(c+t)\Gamma(c-a)}\,,
$$
$$
\tilde{A} = \tilde{A}_t = \tilde{A}(a,c,t) = (a,t)A_t\,,
$$
and
$$
B = B_t = B(a,c,t) = P(a,c,t)-P(a,c,0)\,.
$$

\bigskip

\alku \lem \label{lem:ambm}
Let $f,g$, and $h$ be real valued functions defined on $[0,\infty)$ such
that $f$ is strictly increasing, $f'$ is strictly decreasing,
$0 < g(x) < h(x)$, and $g'(x) \ge h'(x) > 0$ for all $x \in [0,\infty)$.
Let $F(x) = f(g(x))-f(h(x))$. Then

{\rm (1)} $F$ is strictly increasing on $[0,\infty)$. 

In particular,
with notation as in {\ref{pabnot}}, the function $B$ is strictly 
increasing in $t$, so that
$B(a,c,t) \ge 0$ with equality if and only if $t = 0$.

{\rm (2)} $\partial A/\partial c = AB$.
\loppu

\medskip

\proof
(1) By the assumptions,
\begin{eqnarray*}
F'(x) & = & f'(g(x))g'(x) - f'(h(x))h'(x) \\
& > & f'(g(x))g'(x) - f'(g(x))g'(x) = 0\,.
\end{eqnarray*}
We now take $f = \Psi$, $g(x) = c-a+x$, and $h(x) = c+x$. Then by the 
above, $F(x) = \Psi(c-a+x) - \Psi(c+x)$ is strictly increasing on
$[0,\infty)$ so that $F(x) - F(0) \ge 0$ with equality if and only if
$x = 0$. By the definition of $B$, this
means that $B(a,c,t) \ge 0$ with equality if and only if $t = 0$.

(2) By logarithmic differentiation we get
\begin{eqnarray*}
\frac{\partial A/\partial c}{A} & = & 
\frac{\Gamma'(c-a+t)}{\Gamma(c-a+t)}
-\frac{\Gamma'(c+t)}{\Gamma(c+t)}
-\left( \frac{\Gamma'(c-a)}{\Gamma(c-a)}
-\frac{\Gamma'(c)}{\Gamma(c)}\right) \\
& = & \Psi(c-a+t)-\Psi(c+t)-(\Psi(c-a)-\Psi(c)) \\
& = & P(a,c,t) - P(a,c,0) \\
& = & B(a,c,t)\,. \qquad \square
\end{eqnarray*}

\bigskip

\alku\theo \label{quotfdepc}
For $a > 0$ and $x,y \in (0,1)$, the function $f$ defined on 
$(a,\infty)$ by
$$
f(c) = B(a,c-a) \frac{F(a,c-a;c;x)}{F(a,c-a;c;y)}
$$
is strictly decreasing from $(a,\infty)$ onto $(0,\infty)$.
\loppu

\medskip

\proof
First, since
$$
1 \le F(a,c-a;c;x) \le 1+\frac{c-a}{c} F(a,1;1;x),
$$ 
it follows that $F(a;c-a;c;x) \to 1$ as $c \to a^+$. 
Hence
$$ 
f(a+) = \lim_{c \to a+} B(a,c-a) = 
\lim_{c \to a+} \Gamma(c-a) = \infty\,.  
$$

Next, we note that for $n \ge 1,\;\; c \mapsto (c-a,n)/(c,n)$ is increasing by \cite[1.58(32)]{AVV} with limit $1$ as $c \to \infty$.
Hence, using \cite[1.20 (1)]{AVV}, we get
\begin{equation}\label{fatone}
F(a,c-a;c;r) \le F(a,1;1;r) = (1-r)^{-a}.
\end{equation}
Now let $F(c,r)=F(a,c-a;c;r),\; h(r)=(1-r)^{-a}$ and
$$F_n(c,r)=\sum_{k=0}^n \frac{(a,k)(c-a,k)}{(c,k)} \frac{r^k}{k!}\qquad{\rm and} \qquad
h_n(r)=\sum_{k=0}^n (a,k) \frac{r^k}{k!}.$$
Then let $r \in (0,1)$ and $\varepsilon > 0$. Now let $m_0$ be 
such that $h(r)-h_m(r) < \varepsilon$ for all $m>m_0$. Then there 
exists a $c_0$ such that when $c>c_0$ and all $0\le m \le m_0$ we have
$(c-a,m)/(c,m) > 1-\varepsilon$. Thus, for  $c>c_0$ and $p=m_0$ 
we have
\begin{eqnarray*}
F(c,r)&>& F_p(c,r) =\sum_{n=0}^p \frac{(a,n)(c-a,n)}{(c,n)} \frac{r^n}{n!}
>(1-\varepsilon) \sum_{n=0}^p (a,n)\frac{r^n}{n!}\\
&=&(1-\varepsilon) h_p(r) >(1-\varepsilon) (h(r)-\varepsilon).
\end{eqnarray*}  
From this and (\ref{fatone}) we see that $F(a,c-a;c;r) \to (1-r)^{-a}$ 
as $c \to \infty$.
Applying for $z = x,y$ as $c \to \infty$, we see that 
$F(a,c-a;c,x)/F(a,c-a;c;y) \to ((1-x)/(1-y))^{-a}$, which is finite.
As $c \to \infty$, by Stirling's formula \cite[12.33]{WW},
$$
\frac{B(a,c-a)}{\Gamma(a)} = 
\frac{\Gamma(c-a)}{\Gamma(c)} \sim
\left(\left(1-\frac{a}{c}\right)^{c-(1/2)}\right) 
\left(\frac{e}{c-a}\right)^{a} \to (e^{-a})(0) = 0\,.
$$
Hence $f(c) \to 0$, as $c \to \infty$.

Logarithmic differentiation together with the Notation \ref{pabnot} and
Lemma \ref{lem:ambm} (2) yield
\begin{eqnarray*}
\frac{f'(c)}{f(c)} & = & 
\frac{1}{f(c)} \frac{\partial}{\partial c} \left( B(a,c-a)
\frac{F(a,c-a;c;x)}{F(a,c-a;c;y)} \right) \\
 & = & \frac{1}{f(c)} \left( 
 \left(\frac{\partial}{\partial c} \frac{\Gamma(a) \Gamma(c-a)}{\Gamma(c)} \right)
 \frac{F(a,c-a;c;x)}{F(a,c-a;c;y)} \right. \\
 & & \left. + \frac{\Gamma(a) \Gamma(c-a)}{\Gamma(c)} 
 \frac{\partial}{\partial c}
 \frac{F(a,c-a;c;x)}{F(a,c-a;c;y)} \right) \\
 & = & \frac{1}{f(c)} \left( 
 \frac{\Gamma'(c-a) \Gamma(a) \Gamma(c) - \Gamma'(c) \Gamma(a) \Gamma(c-a)}{
 \Gamma(c)^2} \frac{F(a,c-a;c;x)}{F(a,c-a;c;y)} \right. \\
 & & + \frac{\Gamma(a) \Gamma(c-a)}{\Gamma(c)} 
 \frac{1}{F(a,c-a;c;y)^2} \left(
 \left(\frac{\partial}{\partial c} F(a,c-a;c;x)\right) F(a,c-a;c;y) \right. \\
 & & - \left. \left.
 \left(\frac{\partial}{\partial c} F(a,c-a;c;y)\right) F(a,c-a;c;x)
 \right) \right) \\
 & = & \frac{\Gamma'(c-a) \Gamma(a) \Gamma(c) - \Gamma'(c) \Gamma(a) 
 \Gamma(c-a)}{\Gamma(c)^2} \frac{\Gamma(c)}{\Gamma(a)\Gamma(c-a)} \\
 & & + \frac{1}{F(a,c-a;c;x)} 
 \sum_{n=0}^{\infty} \left(\frac{\partial}{\partial c} \tilde{A}_n \right) 
\frac{{x}^{n}}{n!}
 - \frac{1}{F(a,c-a;c;y)} 
 \sum_{n=0}^{\infty} \left(\frac{\partial}{\partial c} \tilde{A}_n \right) 
\frac{{y}^{n}}{n!} 
\\
 & = & \Psi(c-a) - \Psi(c) + \frac{1}{F(a,c-a;c;x)} 
 \sum_{n=0}^{\infty} \tilde{A}_n B_n \frac{{x}^{n}}{n!} \\
 & & - \frac{1}{F(a,c-a;c;y)} 
 \sum_{n=0}^{\infty} \tilde{A}_n B_n \frac{{y}^{n}}{n!}\,.
\end{eqnarray*} 
It follows that
\begin{eqnarray*}
h(c) & = & \frac{1}{B(a,c-a)} F(a,c-a;c;y)^2 f'(c) \\
 & = & F(a,c-a;c;y)F(a,c-a;c;x)(\Psi(c-a) - \Psi(c)) \\
 & & + F(a,c-a;c;y) \sum_{n=0}^{\infty} 
 \frac{\tilde{A}_n B_n}{n!} {x}^{n} 
 - F(a,c-a;c;x) \sum_{n=0}^{\infty} 
 \frac{\tilde{A}_n B_n}{n!} {y}^{n} \\
 & = & (\Psi(c-a) - \Psi(c)) \sum_{n=0}^{\infty} \sum_{m=0}^{\infty}
 \frac{\tilde{A}_n \tilde{A}_m}{n! m!} y^{n} {x}^{m} \\
 & & + \sum_{n=0}^{\infty} \sum_{m=0}^{\infty}
 \frac{\tilde{A}_n \tilde{A}_m B_m}{n! m!} {x}^{m} {y}^{n}
 - \sum_{n=0}^{\infty} \sum_{m=0}^{\infty}
 \frac{\tilde{A}_n \tilde{A}_m B_m}{n! m!} {y}^{m} {x}^{n} \\
 & = & \sum_{n=0}^{\infty} \sum_{m=0}^{\infty}
 \frac{\tilde{A}_n \tilde{A}_m}{n! m!} G_{m,n}(a,c,r) (xy)^{m}\,,
\end{eqnarray*}
where
\begin{eqnarray*}
G_{m,n}(a,c,r) & = & (\Psi(c-a) - \Psi(c))y^{n-m}
+ B_m y^{n - m} - B_m {x}^{n-m} \\
& = & y^{n-m}(\Psi(c-a+m)-\Psi(c+m)) - B_m {x}^{n-m}\,.
\end{eqnarray*}
Since $\Psi$ is strictly increasing, we have
$\Psi(c-a+m)-\Psi(c+m) < 0$ and by Lemma \ref{lem:ambm} (1), $B_m \geq 0$.
Hence $G_{m,n}(a,c,r) < 0$. It follows that $h(c) < 0$ and as
$B(a,c-a) = \Gamma(a)\Gamma(c-a)/\Gamma(c) > 0$ for $0<a<c$,
we get that $f'(c) < 0$ for $c \in (a,\infty)$. $\qquad \square$

\bigskip

\alku \cor \label{mudepc}
For $a>0$ and $r \in (0,1)$ the function $\tilde{f}(c)$ defined on
$(a,\infty)$ by $\tilde{f}(c) = \mu_{a,c}(r)$
is strictly decreasing from $(a,\infty)$ onto $(0, \infty)$
with $\tilde{f}(1)=\mu_a(r)$ if $a<1$.
\loppu

\bigskip

\alku \lem \label{finvmono}
Let $z = f(x,y) = f_x(y) = f_y(x)$ be continuously differentiable for $x$ and 
$y$ in some real intervals. 
Suppose that $(\partial f/ \partial x)(\partial f/ \partial y) > 0$.
Let 
$y = f_x^{-1}(z) = g(x,z)$. Then
$$
\frac{\partial y}{\partial x} = \frac{\partial g}{\partial x} < 0\,.
$$
\loppu

\medskip

\proof
By implicit differentiation 
partial to $x$, we get
$$
0 = \frac{\partial f}{\partial x} +
\frac{\partial f}{\partial y}\frac{\partial g}{\partial x}\,.
$$
Hence
$$
\frac{\partial g}{\partial x} = 
-\frac{\partial f/\partial x}{\partial f/\partial y} < 0\,. \qquad \square
$$

\bigskip

\alku \theo \label{imudpec}
Let 
$a,x > 0$ be fixed. Then the function
$$
g: c \mapsto \mu_{a,c}^{-1}(x)
$$
is strictly decreasing from $(a,\infty)$ onto $(0,1)$ with
$g(1) = \mu_a^{-1}(x)$ if $a < 1$.
\loppu

\medskip

\proof
Denote $r = \mu_{a,c}^{-1}(x) = h(c,x)$. Then
$x = \mu_{a,c}(r) = f(c,r)$.
Now $\partial f/ \partial r < 0$ and by Corollary \ref{mudepc}
$\partial f/\partial c < 0$, so that $\partial g/\partial c < 0$ and the 
monotonicity of $g$ follows from Lemma \ref{finvmono}.

Since $\mu_{a,1} = \mu_a$, we get
$$
x = \mu_a(\mu_a^{-1}(x)) = \mu_{a,1}(\mu_a^{-1}(x))
$$
so that
$$
g(1)=\mu_{a,1}^{-1}(x) = \mu_a^{-1}(x)\,.
$$

We claim that $\lim_{c \to \infty} h(c,x) = 0$. Assume on the contrary that
$\lim_{c \to \infty} h(c,x) = r_0 > 0$. 
Then $h(c,x) > r_0$ for all $c \in (a,\infty)$.
Hence
$$
x = \mu_{a,c}(h(c,x)) < \mu_{a,c}(r_0)\,.
$$
Letting $c \to \infty$, Corollary \ref{mudepc} implies that
$x \le 0$, which is a contradiction.

It remains to show that $h(a+,x)=1$. Suppose that
$h(a+,x) = r_0 \in (0,1)$. 
Then $h(c,x) < r_0$ for all $c \in (a,\infty)$.
Hence
$$
x = \mu_{a,c}(h(c,x)) > \mu_{a,c}(r_0)\,.
$$
Letting $c \to a+$, we get, by Theorem \ref{mudepc}, that
$x = \infty$ which is a contradiction. Thus $h(a+,x) = 1$.
$\qquad \square$

\bigskip

\alku \theo \label{phidepc}
Let $a,r \in (0,1)$ and $K \in (1,\infty)$ be fixed. Then the function
$$
c \mapsto \varphi^{a,c}_{K}(r)
$$
is strictly decreasing from $(a,1]$ onto $[\varphi^{a}_{K}(r), 1)$ and 
the function
$$
c \mapsto \varphi^{a,c}_{1/K}(r)
$$
is strictly increasing from $(a,1]$ onto $(0, \varphi^{a}_{1/K}(r)]$.
\loppu

\medskip

\proof
It is obvious (see \cite[Remark 4.12]{hvv}) that we have 
$\varphi^{a,c}_{K}=\tilde{\mu}_{a,c}^{-1}(\tilde{\mu}_{a,c}(r)/K)$, where
$$\tilde{\mu}_{a,c}(r)=\frac{F(a,c-a;c,r'{\thinspace}^2)}{F(a,c-a;c;r^2}.$$
Denote $s = \varphi^{a,c}_{K}(r)$ and
$$
Q(a,c,r) = F(a,c-a;c;r^2)\,.
$$
By definition,
\be \label{eq:QK}
\frac{Q(a,c,s')}{Q(a,c,s)} = \frac{1}{K} \frac{Q(a,c,r')}{Q(a,c,r)}\,.
\ee
We apply logarithmic differentiation with respect to $c$ to 
(\ref{eq:QK}) and get
$$
\frac{1}{Q(a,c,s')} \left( \frac{\partial Q(a,c,s')}{\partial c} 
- \frac{\partial Q(a,c,s')}{\partial s'} \frac{s}{s'} 
\frac{\partial s}{\partial c} \right)
- \frac{1}{Q(a,c,s)} \left( \frac{\partial Q(a,c,s)}{\partial c} 
+ \frac{\partial Q(a,c,s)}{\partial s} \frac{\partial s}{\partial c} \right)
$$
$$
= \frac{1}{Q(a,c,r')} \frac{\partial Q(a,c,r')}{\partial c}
- \frac{1}{Q(a,c,r)} \frac{\partial Q(a,c,r)}{\partial c},
$$
which is equivalent to
\be \label{eq:QQ1}
\left( \frac{\partial Q(a,c,s')}{\partial s'} \frac{s}{s'} \frac{1}{Q(a,c,s')}
+ \frac{1}{Q(a,c,s)}\frac{\partial Q(a,c,s)}{\partial s} \right) 
\frac{\partial s}{\partial c}
\ee
$$
= (Q_1(a,c,s') - Q_1(a,c,r')) + (Q_1(a,c,r) - Q_1(a,c,s)),
$$
where
$$
Q_1(a,c,x) = \frac{1}{Q(a,c,x)} \frac{\partial Q(a,c,x)}{\partial c}\,.
$$
Then for $0 < a < c \leq 1$ we get that
$$
\frac{\partial Q(a,c,x)}{\partial x} = 
\frac{2 a (c-a)}{c}\ x\ F(a+1,c-a+1;c+1;x^2) > 0
$$
for all $x \in (0,1)$. Hence the coefficient of 
$\partial s/\partial c$ in (\ref{eq:QQ1})
is positive. We turn our attention to the right hand side of (\ref{eq:QQ1}).
By Lemma \ref{lem:ambm} (2) we have that
$$
Q_1(a,c,x) = \frac{\sum_{n=0}^{\infty} \frac{\tilde{A}_n B_n}{n!} r^{2n}}{
\sum_{n=0}^{\infty} \frac{\tilde{A}_n}{n!} r^{2n}}
= \frac{\sum_{n=0}^{\infty} \alpha_n r^{2n}}{
\sum_{n=0}^{\infty} \beta_n r^{2n}}\,,
$$
where $\alpha_n = \tilde{A}_n B_n /n!$ and $\beta_n = \tilde{A}_n / n!$. Then
$$
\frac{\alpha_n}{\beta_n} = B_n = 
\Psi(c) - \Psi(c-a) -(\Psi(c+n) - \Psi(c-a+n))\,,
$$
where, by (\ref{eq:psiser}),
\begin{eqnarray*}
-(\Psi(c+n) - \Psi(c-a+n)) & = & \frac{1}{c+n}
-\sum_{k=1}^{\infty} \frac{c+n}{k(k+c+n)} \\
& & - \frac{1}{c-a+n} 
+ \sum_{k=1}^{\infty} \frac{c-a+n}{k(k+c-a+n)} \\
& = & - \frac{a}{(c+n)(c-a+n)} + \sum_{k=1}^{\infty} 
\frac{-a}{(k+c+n)(k+c-a+n)} \\
& = & -a \sum_{k=0}^{\infty} \frac{1}{(k+c+n)(k+c-a+n)}\,,
\end{eqnarray*}
which is clearly increasing in $n$. Hence $\alpha_n / \beta_n$ is increasing
in $n$ and \cite[Theorem 4.4]{hvv} implies that $Q_1(a,c,x)$ is strictly 
increasing in $x$. Since $K > 1$, it is immediate that $s > r$ and
$r' > s'$ and it follows that the right hand side of (\ref{eq:QQ1}) is
negative. Hence $\partial s/\partial c < 0$, 
which proves the first monotonicity claim.
On the other hand, if $s = \varphi^{a,c}_{1/K} (r)$, then
$s < r$ and $r' < s'$ and the right hand side of (\ref{eq:QQ1}) together
with $\partial s/\partial c$ are positive and the second 
monotonicity claim follows.

It remains to consider the ranges of the functions. The values at
$c = 1$ follow from the fact that for all $k > 0$,
\be \label{eq:phiaphiak}
\phit^{a,1}_{k}(r) = \varphi^a_k(r)\,.
\ee
To show that (\ref{eq:phiaphiak}) holds, we write
$$
\mu_a(\mut_{a,1}^{-1}(t)) = \frac{\pi}{2 \sin (\pi a)}
\mut_{a,1}(\mut_{a,1}^{-1}(t)) = \frac{\pi}{2 \sin (\pi a)} t
$$
and put $t = \mut_{a,1}(r)/k$ to get
$$
\mu_a(\mut_{a,1}^{-1}(\mut_{a,1}(r)/k))
= \frac{\pi}{2 \sin (\pi a)} \frac{\mut_{a,1}(r)}{k}
= \mu_a(r)/k
$$
which implies (\ref{eq:phiaphiak}).

To conclude the proof we need to show that as $c \to a+$,
$\phit^{a,c}_{K}(r) \nearrow 1$ and $\phit^{a,c}_{1/K}(r) \searrow 0$.
We prove the first fact and note that the proof of the second one is
similar. Let $L = \phit^{a,a+}_{K}(r)$. Assume that $L < 1$. 
By the monotonicity in $c,$ it follows that
$\phit^{a,c}_K(r) < L$ for all $c \in (a,1]$.  Hence
$\mut_{a,c}(L) < \mut_{a,c}(r)/K$, so that
$$
\frac{\mut_{a,c}(L)}{\mut_{a,c}(r)} < 1/K\,. 
$$
Letting $c \to a+$, we get $1/K \geq 1$, which is a contradiction, since
$K > 1$. Hence $L = 1$. $\qquad \square$

\bigskip

\alku \theo \label{th:KEB}
For $a, r \in (0,1) ,$ let $f$ and $g$ be functions 
defined on $(a,\infty)$ by

(1) $f(c) = {\K}_{a,c}  - (B/2),$ 

(2) $g(c) = (B/2) - \E_{a,c}$, where $B = B(a,c-a)$. 
Then, both $f$ and $g$ are strictly decreasing, with 
$f(a+) = \log(1/r'), f(\infty) = 0 = g(\infty)$,
$$
g(a+) = \frac{1}{2} \left(\sum_{n= 1}^{\infty} \frac{r^{2n} }{a+n-1} 
\right) - \log(1/r') \,. 
$$
\loppu

\medskip

\proof 
The assertion $f(\infty) = 0= g(\infty)$ follows immediately 
from Stirling's formula, as in the proof of Theorem \ref{quotfdepc} (cf. \cite[1.49.]{AVV}). Next, the coefficient of $r^{2n}$ in the Maclaurin series 
of $f(c)$ is
$$
f_n(c) = \Gamma(a+n)\Gamma(c-a+n)/(2 (n!)\Gamma(c+n))\,, 
$$ 
so that
$f_n'(c)/f_n(c) = \Psi(c-a+n) - \Psi(c+n) < 0$, since $\Psi$ is strictly 
increasing. Similarly, it can be shown that an analogous assertion holds 
for $g(c)$, thus proving the monotonicity of these functions. Finally,
$$
f(a+) = \sum_{n=1}^{\infty}{ \frac{r^{2n}}{2n} } = \log(1/r') \,,
$$
and
\begin{eqnarray*}
g(a+) & = & 
\frac{1}{2} \sum_{n=1}^{\infty} {  \frac{(1-a)r^{2n}}{n(a+n-1)}} \\
 & = & - \frac{1}{2} \sum_{n=1}^{\infty}\left(\frac{ r^{2n}}{n} - 
\frac{r^{2n}}{a+n-1} \right) \\
 & = & \frac{1}{2} \sum_{n= 1}^{\infty}\left( \frac{r^{2n}}{a+n-1} \right) 
- \log(1/r')\,. \qquad \square
\end{eqnarray*}

Finally, we make some conjectures regarding the behavior of the Legendre
$\M$-function combined with other functions. Such problems seem to be quite 
difficult, and apart from the functions in Theorem \ref{Mfunctions} and 
immediate consequences, we are not aware of any results in this direction.
In particular, solving any one of the following problems immediately yields
several interesting functional inequalities generalizing those stated in
\cite[1.14, 1.15]{AQVV}.   

\alku \conj \label{ops}
Based on experimental evidence, we make the following conjectures.\\
(1)  Let $0<a<c<1$. Then the function $f(r)=\sqrt{r}/\M(r^2)$ is 
strictly increasing from $(0,1)$ onto $(0,B)$, and $g(r)=\sqrt{r'}/\M(r^2)$
is strictly decreasing from  $(0,1)$ onto $(0,B)$.\\
(2) Let $0<a<c<1,\; K >1,$ and $s=\varphi_K^{a,c}(r)$. Then the function\\
\hspace*{5mm}{\it i)} $f_1(r)=(s \M(r^2))/(r \M(s^2))$ is decreasing
from $(0,1)$ onto $(1,\infty)$.\\
\hspace*{5mm}{\it ii)} $f_2(r)=(s' \M(r^2))/(r' \M(s^2))$ is decreasing
from $(0,1)$ onto $(0,1)$.\\
\hspace*{5mm}{\it iii)} $f_3(r)=(\K(r) \M(r^2))/(\K(s) \M(s^2))$ is decreasing
from $(0,1)$ onto $(1/K,1)$.\\  
\hspace*{5mm}{\it iv)} $f_4(r)=(\K'(r) \M(r^2))/(\K'(s) \M(s^2))$ is decreasing from $(0,1)$ onto 
\hspace*{5mm}$(1,K)$.\\  
\loppu

\bigskip

\end{document}